\documentclass[11pt]{article}

\usepackage[T1]{fontenc}
\usepackage{graphicx}
\usepackage{amsmath,amssymb,bm}
\usepackage{hyperref}
\usepackage[nameinlink,capitalize,noabbrev]{cleveref}
\usepackage{color}
\usepackage{blindtext}
\usepackage{enumitem}
\usepackage[ruled,linesnumbered]{algorithm2e}
\usepackage{amsthm}

\newtheorem{definition}{Definition}[section]

\theoremstyle{remark}

\begin{document}

\title{Devil's terraces: determining the organization of resonance tongues in a periodically forced dynamical system}

\author{John Bailie, Priya Subramanian, Bernd Krauskopf \\ Department of Mathematics, University of Auckland, Auckland, New Zealand}
\date{}
\maketitle

\begin{abstract}
\noindent
In periodically forced dynamical systems, resonance tongues are open regions of a parameter plane in which the dynamics on an invariant torus locks to a stable periodic orbit. While individual resonance tongues are well understood, the principles governing their global arrangement remain largely unexplored. We develop a topological framework, grounded in applied topology and Morse theory, whose central object is the two-dimensional \emph{resonance surface}, defined as the graph of the rotation number~$\rho$ over a parameter plane. Within this framework, resonance tongues appear as terraces of the resonance surface at rational values of~$\rho$, and their global arrangement is determined by the singularities of this surface. Resolving the resonance surface requires the accurate computation of~$\rho$, and we present an algorithm that does so efficiently and at high resolution. As a specific example, we examine a periodically forced model of vertical mixing in the North Atlantic, a process relevant to the Atlantic Meridional Overturning Circulation, and study how its resonance surface changes under variation of a third parameter. We identify six distinct resonance-tongue arrangements and show that the \emph{resonance transitions} between them are due to changes in the number and type of singularities on the boundary of the resonance surface.
\end{abstract}
\maketitle

\section{Introduction}
\label{sec:intro}
Periodically forced models arise across a wide range of physical systems, where they exhibit multi-frequency dynamics and resonance phenomena, such as mode locking in biophysical neuron models~\cite{farokhniaee2017mode,lee2006bifurcation,wang2014influence}, synchronization in injection laser systems~\cite{terrien2023merging,simonet1994locking}, phase locking of delay-oscillator models of the El~Ni\~no--Southern Oscillation (ENSO) to the seasonal cycle~\cite{keane2018chenciner,keane2017climate,tziperman1995irregularity,krauskopf2014bifurcation}, and resonant spatial patterns in the Belousov--Zhabotinsky reaction under time-periodic perturbation~\cite{lin2004resonance,marts2007period}.

From a dynamical systems perspective, resonance in a periodically forced differential equation is the result of an interaction between the external forcing and a periodic orbit of the unforced system. For sufficiently weak forcing, an unforced periodic orbit generically becomes a normally hyperbolic two-dimensional invariant torus that inherits its stability~\cite{wiggins2003introduction}. The rotation number~$\rho$ on this torus distinguishes two qualitatively different dynamical regimes~\cite{kuznetsov1998elements}: if $\rho = p/q \in \mathbb{Q}$, the system is in $p{:}q$ resonance, meaning that trajectories on the torus converge to a stable period-$q$ orbit that forms a $p{:}q$ torus knot, while an unstable period-$q$ orbit of the same type is also present; if~$\rho$ is irrational, the dynamics are quasi-periodic and every trajectory is dense on the torus~\cite{shil2001methods}. In a two-parameter plane, rational values of~$\rho$ occur in open regions known as \emph{resonance} (or \emph{Arnold}) \emph{tongues}, each bounded by saddle-node bifurcation curves along which the stable and unstable $p{:}q$ periodic orbits coalesce; irrational values of~$\rho$ are attained along curves separating these resonance tongues~\cite{arnold1965small,kuznetsov1998elements}.

The theory of individual resonance tongues is well developed. It began with Arnold's analysis of a family of circle maps~\cite{arnold1965small}, which provided the prototype for resonance tongues emerging from a zero-coupling line at rational rotation numbers~\cite{broer2000resonance}. Subsequent work developed equivariant normal forms for analyzing strong resonances~\cite{takens2001forced,arnold1977loss} and applied singularity theory to classify tongue-boundary geometry~\cite{broer2003geometry}. Studies of how resonance tongue geometry deforms under parameter variation have shown that resonance regions can pinch off into closed, banana-shaped regions~\cite{peckham1995bananas} or develop intricate `Arnold flame' structures~\cite{peckham2002lighting}; continuation-based methods later made it possible to represent the resonance tongues themselves as surfaces over parameter space~\cite{schilder2007computing}.

The global arrangement of resonance tongues has received far less attention. Terrien et al.~\cite{terrien2023merging}, studying a delay differential equation model of a pulsing excitable micro-laser, showed that extrema of the rotation number~$\rho$ along a torus bifurcation curve organize nearby resonance tongues and mark the points at which tongues merge or disconnect as a third parameter is varied. A similar organizing role of extrema of~$\rho$ has been observed in delay differential equation models of the El~Ni\~no--Southern Oscillation~\cite{keane2018chenciner, bolduc2025seasonal}. Interactions between individual resonance tongues have also been reported away from torus bifurcation, including their coalescence under one-parameter variation~\cite{habib2018isolated, detroux2018experimental, kuether2015nonlinear} and their appearance as isolated components (isolas) in a two-parameter plane~\cite{marchionne2018synchronisation}. However, extrema of~$\rho$ were not identified as the organizing mechanism in these studies.

In this paper, we develop a framework, grounded in applied topology and piecewise-linear variants of Morse theory~\cite{edelsbrunner2010computational, carr2003computing}, for describing the global arrangement of resonance tongues in a two-parameter plane. Our central object of study is the \emph{resonance surface}~$\mathcal{S}$, defined as the graph of the rotation number~$\rho$ over the relevant region of the two-parameter plane where one finds a normally hyperbolic invariant torus. In this formulation, resonance tongues correspond precisely to the level sets of~$\rho$ at rational values, and the question of how resonance tongues are arranged becomes one of understanding the topology of these level sets. It is well known that crossing resonance tongues along a one-dimensional path in the parameter plane gives a devil's staircase~\cite{wiggins2003introduction}. In our context, we are hence interested in the topological nature of the ``devil's terraces'' encoded by the resonance surface at rational values of~$\rho$.

\begin{figure}[t!]
  \centering
  \includegraphics{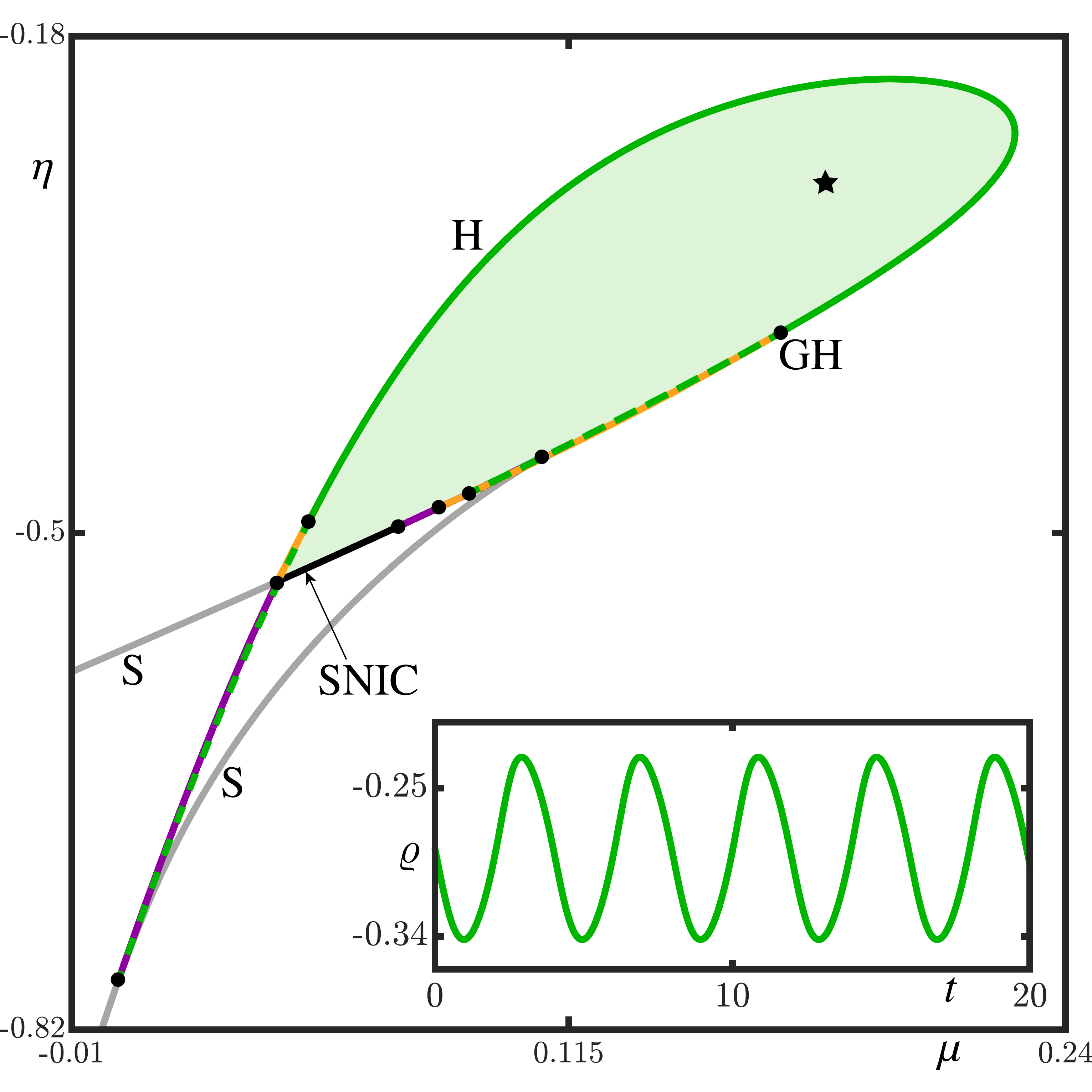}
\caption{\label{fig:bd_nof}Bifurcation diagram of system~\eqref{eq:forced_ODE_second} in the $(\mu,\eta)$-plane at $c=0$, for $\kappa_1=0.1$, $\kappa_2=1.0$, and $\varepsilon=0.05$. The shaded (green) region indicates where the global attractor is a single stable periodic solution; this region is bounded by a curve~$\mathrm{H}$ of Hopf bifurcations (supercritical: solid green; subcritical: dotted green) and by a segment~$\mathrm{SNIC}$ (black) of saddle-node bifurcations on an invariant circle along the larger saddle-node curve~$\mathrm{S}$ (gray). Also shown is the generalized Hopf point~$\mathrm{GH}$ at which~$\mathrm{H}$ changes criticality. The further colored curves and black codimension-two points are summarized in~\cite{bailie2024bifurcation}. The inset shows the density oscillation~$\varrho(t)$ at the parameter point marked by the star.}
\end{figure}

For smooth surfaces, Morse theory~\cite{milnor1963morse, audin2014morse, matsumoto2002introduction, cerf1970stratification} provides the classical framework for relating the number and position of critical points to changes in level-set topology. The rotation number~$\rho$, however, is not a Morse function: it is locally constant on open regions where it takes rational values, and it varies monotonically between them. The classical theory therefore does not apply directly, and the notion of a critical point must be generalized. In particular, the critical set of~$\rho$ need not consist of isolated points but may extend along entire resonance tongues. To examine~$\mathcal{S}$, we develop theoretical tools to identify and classify the critical values of~$\rho$ --- referred to here as \emph{singularities} --- and use them to decompose~$\mathcal{S}$ into specific equivalence classes of resonance tongues. This decomposition provides a combinatorial description of structurally stable resonance tongue configurations. It is also a basis for classifying \emph{resonance transitions} as codimension-one events at which the qualitative arrangement of resonance tongues in a parameter plane changes under variation of a third parameter.

This framework is motivated by and applied to a periodically forced version of a planar ordinary differential equation model due to Welander~\cite{welander1986thermohaline, cessi1994simple, bailie2025detailed, bailie2024bifurcation}, which describes vertical mixing in the North Atlantic and is relevant to the Atlantic Meridional Overturning Circulation. This model, which we refer to as the periodically forced Welander model, depends on three parameters: the virtual salinity flux~$\mu$, a density threshold~$\eta$, and the forcing amplitude~$c$. For a fixed value of~$\eta$, the system admits a two-dimensional invariant torus over a half-disk-shaped region of the $(\mu,c)$-plane, bounded above by a torus bifurcation curve and below by the segment where $c = 0$. Throughout the range of~$\eta$ considered here, this torus remains normally hyperbolic across the entire region, which makes the model an ideal setting for a case study into how singularities of~$\rho$ organize the structurally stable arrangement of resonance tongues, as captured by the level-set topology of the resonance surface~$\mathcal{S}$. These singularities occur both on the boundary of the half-disk and in its interior. For fixed~$\eta$, the number and type of boundary singularities organize the resonance tongues in the $(\mu,c)$-plane into six qualitatively distinct, structurally stable configurations, denoted cases~$\mathrm{I}$--$\mathrm{VI}$. We focus here on the resonance transitions between them, which change the boundary singularity structure. However, we will see that cases~$\mathrm{III}$--$\mathrm{V}$ additionally contain several interior singularities. The dependence of these interior singularities on~$\eta$, together with the associated interior resonance transitions, will be left for future work.

Resolving the resonance surface~$\mathcal{S}$ accurately enough to detect its singularities is a significant computational challenge. Resonance tongues with rotation number $\rho = p/q$ become extremely narrow as~$q$ increases and, for the periodically forced Welander model, continuation of $p{:}q$ locked periodic solutions becomes impractical already for $q \geq 14$. Moreover, neither the associated circle map nor its lift is available in closed form. We overcome these difficulties by developing an algorithm that computes~$\rho$ on a fine grid of fixed parameter values directly from the stroboscopic map, thereby resolving tongues at large~$q$ (without continuation of periodic orbits). Along the torus bifurcation curve,~$\rho$ is determined from the Floquet multipliers of the bifurcating periodic orbit, and along the zero-forcing segment it is obtained from the period of the corresponding periodic orbit of the unforced system. The saddle-node bifurcation curves bounding resonance tongues with $\rho = p/q$ and $q \leq 13$ are computed by continuation of the corresponding periodic orbits; higher-order tongues with $14 \le q \leq 10^4$ are resolved instead as level sets of~$\rho$.

The paper is organized as follows. Section~\ref{sec:welander} introduces the periodically forced Welander model and its parameter setting. Section~\ref{section:rot_number_comp} presents the rotation-number algorithm, which computes~$\rho$ from stroboscopic iterates and provides the numerical basis for resolving the resonance surface. Section~\ref{section:casesIandII} presents three motivating examples that illustrate the main geometric features of resonance tongues and the resonance surface. Section~\ref{section:topo_framework} develops our topological framework, which we use to classify the singularities of~$\mathcal{S}$, describe the local arrangement of resonance tongues around them, and define the resonance transitions. Section~\ref{section:res_bd_unfold} presents these resonance transitions by comparing the local resonance structure on either side of each. Section~\ref{section:res_bd} tracks the singularities on the boundary and at corner points in the three-parameter $(\mu,\eta,c)$-space, to locate the transitions that partition the $\eta$-axis into cases~$\mathrm{I}$--$\mathrm{VI}$. Section~\ref{section:cases_organisation} completes the catalogue of resonance-tongue arrangements by describing additional cases~$\mathrm{IV}$--$\mathrm{VI}$. Finally, we draw conclusions and point out directions for future research in Section~\ref{section:conclusion}.

\section{The Welander model}
\label{sec:welander}
The Welander model~\cite{welander1986thermohaline, cessi1996convective} is an autonomous box model for the Atlantic Meridional Overturning Circulation that idealizes a vertical water column in the North Atlantic as two stacked, well-mixed boxes representing the surface ocean and the deep ocean. It describes the evolution of temperature and salinity in the surface box, while the deep box is held at a fixed reference state. We consider a sinusoidally forced extension of the model that incorporates seasonal variability in North Atlantic salinity~\cite{yashayaev2016recurrent, holliday2015multidecadal}. The non-dimensional periodically forced Welander model is given by
\begin{equation}
\label{eq:forced_ODE_second}
\left\{
\begin{aligned}
    \frac{dx}{dt} &= 1 - x - \mathbf{K}_\varepsilon(x,y)\,x, \\[2pt]
    \frac{dy}{dt} &= \mu\bigl(1 + c\cos(2\pi\theta)\bigr) - \mathbf{K}_\varepsilon(x,y)\,y, \\[2pt]
    \frac{d\theta}{dt} &= 1,
\end{aligned}
\right.
\end{equation}
with
\begin{equation*}
\mathbf{K}_\varepsilon(x,y)
=
\kappa_1 + \frac{1}{2}(\kappa_2 - \kappa_1)
\left(1 + \tanh\!\left(\frac{y - x - \eta}{\varepsilon}\right)\right),
\end{equation*}
where $x$ and $y$ are the non-dimensional temperature and salinity in the surface box, and $\theta$ is the phase of the periodic forcing. The non-dimensional density contrast is $\varrho = y - x$, and the forcing has period one, corresponding to variability over one year. The function~$\mathbf{K}_\varepsilon(x,y)$ is a smooth nonlinear ramp between a weak-mixing rate~$\kappa_1$ and a strong-mixing rate~$\kappa_2$, with $0 < \kappa_1 < \kappa_2$; it parametrizes the exchange between the surface and deep boxes. Following~\cite{bailie2024bifurcation, bailie2025detailed, cessi1996convective}, we fix throughout $\kappa_1 = 0.1$, $\kappa_2 = 1.0$, and $\varepsilon = 0.05$. Note that we have written~\eqref{eq:forced_ODE_second} as an autonomous system with phase space~$\mathbb{R}^2\times\mathbb{S}^1$, where~$\mathbb{S}^1 = \mathbb{R}/\mathbb{Z}$.

Figure~\ref{fig:bd_nof} shows the bifurcation diagram of system~\eqref{eq:forced_ODE_second} in the $(\mu, \eta)$-plane for the unforced case $c = 0$. The green-shaded region indicates parameter values for which a stable periodic orbit is the global attractor; the inset shows a corresponding example time series~$\varrho(t)$ over five periods. This region is bounded primarily by the Hopf bifurcation curve~$\mathrm{H}$ and by the segment~$\mathrm{SNIC}$ of saddle-node bifurcations on an invariant circle~\cite{kuznetsov1998elements}, which forms part of the larger saddle-node curve~$\mathrm{S}$. Several further bifurcation curves, together with the codimension-two points at which they meet (black dots), bound narrower regions of the diagram; they are not discussed here, and we refer to Refs.~\cite{bailie2024bifurcation, bailie2025detailed} for a detailed account.

We restrict attention to values of~$\eta$ for which the horizontal line in the $(\mu, \eta)$-plane intersects the curve~$\mathrm{H}$ above the generalized Hopf bifurcation point~$\mathrm{GH}$, where the Hopf bifurcation along~$\mathrm{H}$ is supercritical. For any such fixed~$\eta$, the line intersects~$\mathrm{H}$ at two points; these two Hopf bifurcation points are the endpoints, at $c = 0$, of a torus bifurcation curve in the $(\mu, c)$-plane. This curve bounds a half-disk-shaped region in which the stable periodic orbit of the unforced system perturbs to a normally hyperbolic invariant torus~$\mathcal{T}$ of system~\eqref{eq:forced_ODE_second}.

\begin{figure*}[t!]
  \centering
  \includegraphics{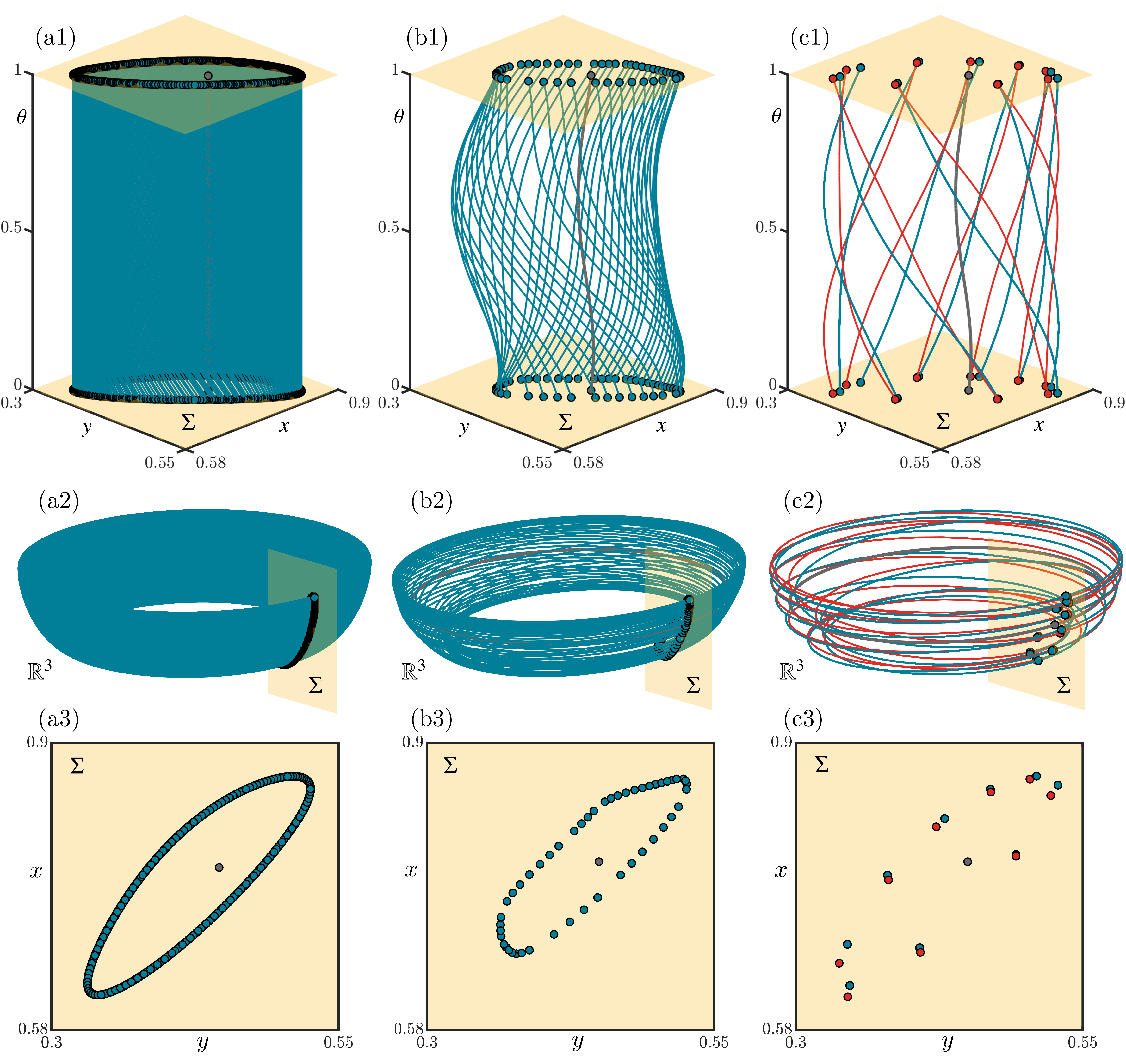}
    \caption{\label{fig:intro_torus}
Dynamics of system~\eqref{eq:forced_ODE_second} on stable invariant tori, shown as (top) trajectories in the cylindrical phase space~$\mathbb{R}^2\times\mathbb{S}^1$, (middle) an embedding in~$\mathbb{R}^3$, and (bottom) iterates of the stroboscopic map on the section~$\Sigma$. In each representation, the gray curve is the unstable period-one orbit. Panels~(a1)--(a3), for $(\mu,c,\eta)=(0.14,\, 0.0,\, -0.286)$, show the invariant torus at zero forcing. Panels~(b1)--(b3), for $(\mu,c,\eta)=(0.14,\, 0.9,\, -0.286)$, show a quasi-periodic torus. Panels~(c1)--(c3), for $(\mu,c,\eta) = (0.14,\,0.34,\,-0.286)$, show a resonant torus with a stable $2{:}9$ periodic orbit, together with the accompanying unstable $2{:}9$ orbit shown in red.}
\end{figure*}

Figure~\ref{fig:intro_torus} illustrates three dynamical regimes on~$\mathcal{T}$, each shown in three complementary views: the orbit in the cylindrical phase space $\mathbb{R}^2 \times \mathbb{S}^1$, the orbit embedded in~$\mathbb{R}^3$ via toroidal coordinates, and the stroboscopic iterates on a Poincar\'e section~$\Sigma$ defined by~$\theta = 0$. Panels~(a1)--(a3) show the unforced case $c = 0$, for which~$\mathcal{T}$ is the trivial suspension of the unforced periodic orbit under the flow of system~\eqref{eq:forced_ODE_second}. Panels~(b1)--(b3) show a quasi-periodic orbit (or locked orbit of very high period) densely winding around~$\mathcal{T}$, whose stroboscopic iterates begin to fill out an invariant circle~$\mathcal{I}$ in~$\Sigma$. Panels~(c1)--(c3) show an attracting and a saddle $2{:}9$ periodic orbit, whose stroboscopic iterates form period-$9$ orbits on~$\mathcal{I}$, corresponding to rotation number~$\rho = 2/9$. In what follows, we focus on the organization across the $(\mu, c)$-plane of the resonance tongues that determine the dynamics on~$\mathcal{T}$. 

\section{Computing the rotation number $\rho$}
\label{section:rot_number_comp}
We now present an algorithm that computes the rotation number~$\rho$ on a normally hyperbolic invariant torus~$\mathcal{T}$, formulated for system~\eqref{eq:forced_ODE_second}. Since the forcing has period one, the dynamics on~$\mathcal{T}$ is sampled by the stroboscopic map on the Poincar\'e section
\begin{equation}
    \Sigma = \{(\mathbf{x},\theta) \in \mathbb{R}^2 \times \mathbb{S}^1 \mid \theta = 0\},
\end{equation}
where $\mathbf{x}=(x,y) \in \mathbb{R}^2$. Let
\begin{equation}
    \varphi_t : \mathbb{R}^2 \times \mathbb{S}^1 \to \mathbb{R}^2 \times \mathbb{S}^1
\end{equation}
denote the flow of~\eqref{eq:forced_ODE_second}. The associated stroboscopic map is
\begin{equation}
    \mathcal{P} : \Sigma \to \Sigma, \qquad \mathcal{P}(\mathbf{x}) = \varphi_{1}(\mathbf{x},0),
\end{equation}
which advances each point of~$\Sigma$ by one forcing period. The torus~$\mathcal{T}$ intersects~$\Sigma$ in an invariant circle
\[
    \mathcal{I} = \mathcal{T} \cap \Sigma,
\]
so determaining~$\rho$ reduces to computing the average angular displacement of the iterates of~$\mathcal{P}$ on~$\mathcal{I}$. The construction extends to higher-dimensional periodically forced systems after projecting~$\mathcal{I}$ onto a plane on which its angular ordering is well defined.

\begin{figure}[ht!]
    \centering      \includegraphics{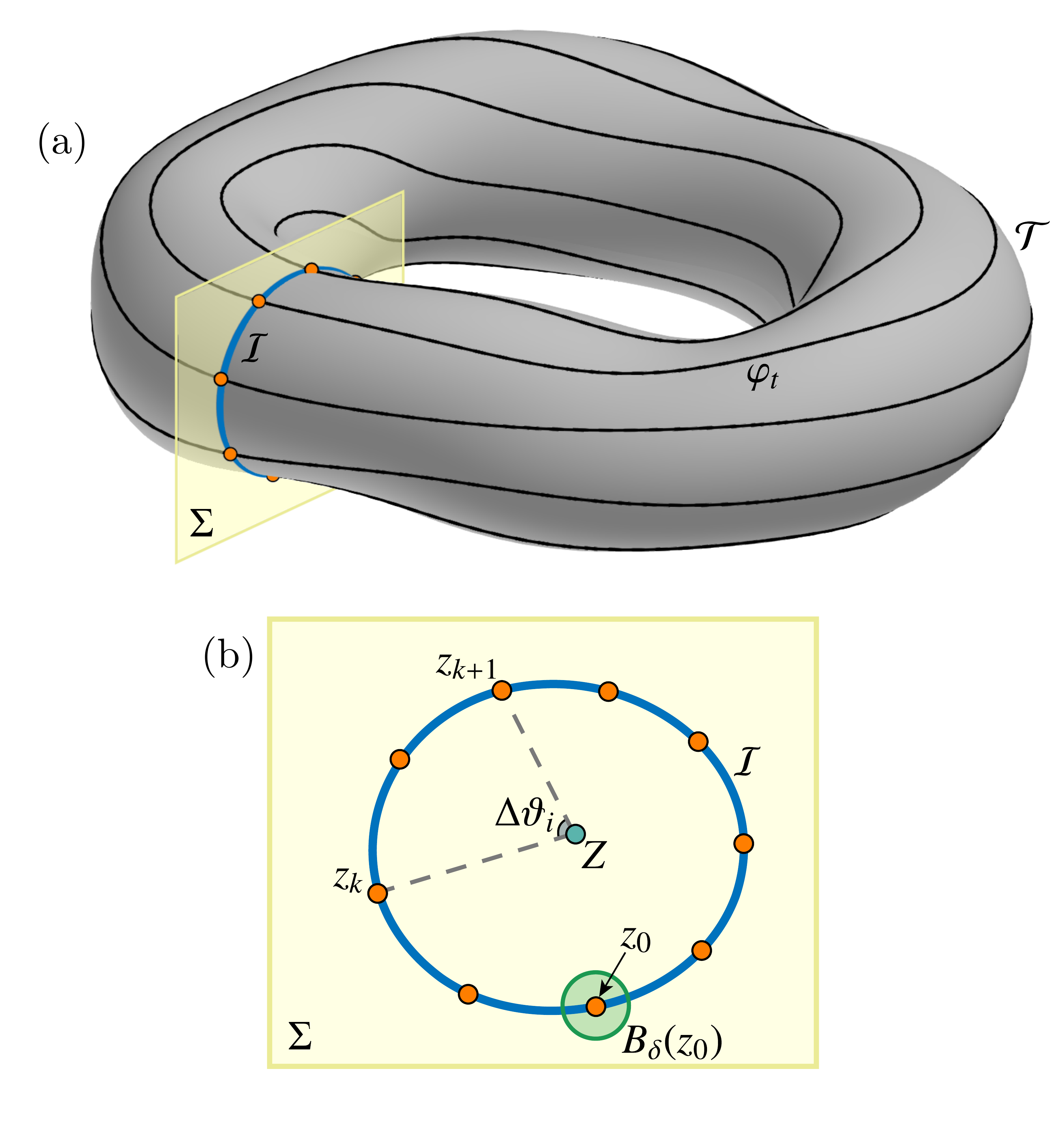}
\caption{\label{fig:torus_algo_example}
Computation of the rotation number~$\rho$ for a $2{:}9$ periodic orbit. Panel~(a) shows a stable invariant torus~$\mathcal T$ (gray) and its intersection~$\mathcal{I}$ (blue) with the stroboscopic Poincar\'e section~$\Sigma$ (yellow plane). Panel~(b) shows successive iterates~$z_k$ (orange dots) of the stroboscopic map, the reference point~$Z$ (green dot), and the angular displacement~$\Delta\vartheta_k$ between consecutive iterates~$z_{k}$ and~$z_{k+1}$. The open ball~$B_\delta(z_0)$ (green disk) is also shown.}
\end{figure}

The algorithm consists of three steps and is illustrated by the schematic in Figure~\ref{fig:torus_algo_example} for a $2{:}9$ locked orbit. Panel~(a) shows the corresponding periodic orbit on the torus~$\mathcal{T}$ embedded in~$\mathbb{R}^3$, and panel~(b) shows the period-$9$ orbit of the stroboscopic map~$\mathcal{P}$ on the invariant circle~$\mathcal{I}$.

In the first step, we choose an initial condition~$\mathbf{x}_0 \in \Sigma$ in the basin of attraction of~$\mathcal{T}$ and integrate system~\eqref{eq:forced_ODE_second} for an integer transient time~$T_{\mathrm{tr}} \in \mathbb{N}$. Since $T_{\mathrm{tr}}$ is an integer number of forcing periods, the endpoint
\[
    z_0 := \varphi_{T_{\mathrm{tr}}}(\mathbf{x}_0)
\]
lies in~$\Sigma$ and, after the transient, is sufficiently close to~$\mathcal{I}$.

In the second step, starting from~$z_0 \in \mathbb{R}^2$, we iterate the stroboscopic map to generate
\begin{equation}
    z_k := \mathcal{P}^k(z_0), \qquad k = 0, \ldots, N,
\end{equation}
where $N \geq 1$ is the first index for which $z_N$ lies in the open ball~$B_\delta(z_0)$ of radius~$\delta > 0$ centered at~$z_0$, as shown in Figure~\ref{fig:torus_algo_example}(b). In the phase-locked case, the return is exact up to numerical accuracy and~$N$ is the period of the orbit of~$\mathcal{P}$. In the quasi-periodic case, no exact return occurs, but recurrence on~$\mathcal{I}$ guarantees arbitrarily close near-returns to~$z_0$. For practical reasons, we impose an iteration cutoff~$N_{\max}$ and return nothing if no near-return is found before this cutoff is reached.

In the third step, we accumulate the angular displacement over the sequence~$\{z_0, \ldots, z_N\}$. We choose a reference point~$Z$ in the interior of the region enclosed by~$\mathcal{I}$; in practice, and as shown in Figure~\ref{fig:torus_algo_example}(b), we take~$Z$ to be the centroid
\begin{equation}
    Z = \dfrac{1}{N}\sum_{k=0}^{N-1} z_k
\end{equation}
of the iterates~$\{z_0, \ldots, z_{N-1}\}$. For each pair of consecutive iterates, let
\begin{equation}
    \Delta\vartheta_k = \angle(z_k,Z,z_{k+1}), \qquad k = 0, \ldots, N-1,
\end{equation}
denote the signed angle from~$z_k$ to~$z_{k+1}$ about~$Z$. The rotation number is then approximated by the Birkhoff average~\cite{das2017quantitative}
\begin{equation}
    \rho = \frac{1}{2\pi N}\sum_{k=0}^{N-1} \Delta\vartheta_k.
    \label{eq:rot_num_form}
\end{equation}
Note that~\eqref{eq:rot_num_form} is exact in the phase-locked case, when $z_N = z_0$, and provides a good finite-time approximation in the quasi-periodic case, with accuracy controlled by~$\delta$ and~$N_{\max}$. The full procedure is summarized in Algorithm~\ref{algo}.

\begin{algorithm}[H]
\label{algo}
\caption{Computing the rotation number $\rho$}
\KwIn{initial condition $\mathbf{x}_0$, transient time $T_{\mathrm{tr}}\in\mathbb{N}$, tolerance $\delta>0$, iteration cutoff $N_{\max}$}
\KwOut{rotation number $\rho\in\mathbb{R}$, or no value}
\BlankLine
\textbf{Step 1: Transient integration.}\\
Integrate system~\eqref{eq:forced_ODE_second} from~$\mathbf{x}_0$ over $[0,T_{\mathrm{tr}}]$ so that the trajectory settles onto~$\mathcal{T}$. Set $z_0 := \varphi_{T_{\mathrm{tr}}}(\mathbf{x}_0)$.
\BlankLine
\textbf{Step 2: Stroboscopic iteration and return index.}\\
\For{$k = 1$ \KwTo $N_{\max}$}{
    $z_k \leftarrow \mathcal{P}(z_{k-1})$\;
    \If{$z_k \in B_\delta(z_0)$}{
        \textbf{break}\;
    }
}
$N \leftarrow k$\;
\If{$N \leq 2$ \textbf{ or } $N=N_{\max}$}{
    \KwRet{no value}
}
\BlankLine
\textbf{Step 3: Computation of~$\rho$.}\\
$Z \leftarrow \dfrac{1}{N}\sum_{k=0}^{N-1} z_k$ \tcp*{reference point inside the region enclosed by~$\mathcal{I}$}
$\rho \leftarrow 0$\;
\For{$k=0$ \KwTo $N-1$}{
    $\Delta\vartheta_k \leftarrow \angle(z_k,Z,z_{k+1})$\;
    $\rho \leftarrow \rho + \Delta\vartheta_k$\;
}
$\rho \leftarrow \rho/(2\pi N)$\;
\BlankLine
\KwRet{$\rho$}
\end{algorithm}

The tolerance~$\delta$ controls how close an iterate of~$\mathcal{P}$ must come to~$z_0$ before it is accepted as a return, and the cutoff~$N_{\max}$ sets the largest return time that can be resolved. In the phase-locked case, the orbit closes after~$N$ iterates and Algorithm~\ref{algo} returns $\rho = p/N$ up to integration error, provided~$\delta$ is sufficiently small and~$N_{\max}$ sufficiently large.

Algorithm~\ref{algo} returns no value (\texttt{nothing} in the Julia implementation) in two situations. The first is when the cutoff $N = N_{\max}$ is reached before the trajectory returns to~$B_\delta(z_0)$, indicating that~$N_{\max}$ is too small or~$\delta$ is too tight. The second is when the trajectory has converged to a fixed point or to a period-$2$ orbit of~$\mathcal{P}$: in the former case the iteration does not converge to an invariant circle~$\mathcal{I}$, while in the latter the two iterates and their centroid~$Z$ are collinear, so the angular displacement about~$Z$ degenerates.

When investigating resonance tongues in system~\eqref{eq:forced_ODE_second}, we compute~$\rho$ on a uniform $800 \times 800$ grid in the $(\mu,c)$-plane. Throughout this region, the invariant torus~$\mathcal{T}$ is the global attractor whenever it exists, so any fixed initial condition in its basin of attraction converges to~$\mathcal{T}$. We use the accuracy settings
\[
    T_{\mathrm{tr}} = 10^3, \qquad
    \delta = 10^{-4}, \qquad
    N_{\max} = 10^4,
\]
which give consistent rotation-number computations across the entire region, including for grid points close to the torus bifurcation curve~$\mathrm{T}$, where the torus is small and normal attraction is weak.

\section{Cases~$\mathrm{I}$--$\mathrm{III}$ as motivating examples}
\label{section:casesIandII}
We now consider the first three resonance-tongue arrangements of system~\eqref{eq:forced_ODE_second}, encountered as~$\eta$ is decreased and denoted cases~$\mathrm{I}$, $\mathrm{II}$, and~$\mathrm{III}$. Each is defined by the number and type of singularities of the rotation number~$\rho$ on the torus bifurcation curve~$\mathrm{T}$ and on the zero-forcing segment~$\mathrm{T}_0$. These cases motivate the topological framework developed in Section~\ref{section:topo_framework}.

\begin{figure*}[ht!]
  \centering
  \includegraphics{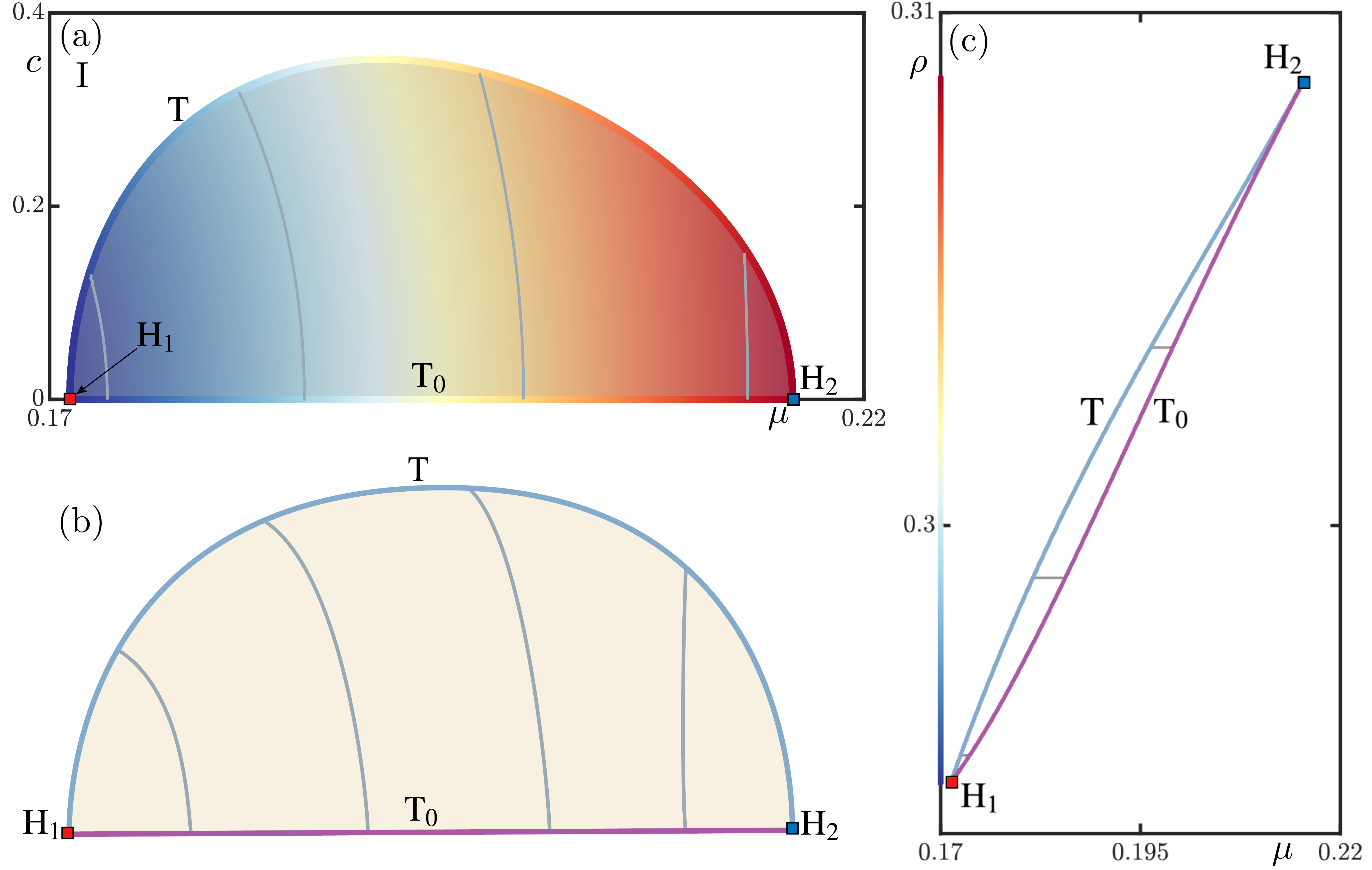}
\caption{\label{bd:fig1} The resonance structure for case~$\mathrm{I}$ at $\eta=-0.214$ on the half-disk~$D$, with the curves~$\mathrm{T}$ and~$\mathrm{T}_0$, which intersect at the corner points~$\mathrm{H}_1$ (red square) and~$\mathrm{H}_2$ (blue square). Panel~(a) shows the resonance diagram in the $(\mu,c)$-plane, with the half-disk colored by the rotation number~$\rho$ from low values in dark blue to high values in dark red; selected resonance tongues are indicated in gray. Panel~(b) shows the corresponding decomposition diagram, with the single interior class shaded pale yellow. Panel~(c) shows the boundary diagram, whose vertical axis is colored by the range of~$\rho$.}
\end{figure*}

We begin with case~$\mathrm{I}$ to introduce the representations and notation that are used throughout. Figure~\ref{bd:fig1} shows the arrangement of resonance tongues in three complementary ways:
\begin{enumerate}[label=(\alph*)]
    \item a \emph{resonance diagram} in the half-disk-shaped region~$D_\eta$ of the $(\mu,c)$-plane, over which a family of stable normally hyperbolic invariant tori exists, so that the rotation number~$\rho$ is defined. The boundary~$B_\eta$ of~$D_\eta$ consists of the zero-forcing segment~$\mathrm{T}_0$ at~$c=0$ and the torus bifurcation curve~$\mathrm{T}$ for~$c>0$, which meet at the Hopf bifurcations~$\mathrm{H}_1$ and~$\mathrm{H}_2$, referred to as the \emph{corner points}. Where no confusion is possible, we suppress the subscript~$\eta$ and write~$D$ and~$B$. The region~$D$ and its boundary~$B$ are colored by~$\rho$, and selected resonance tongues are also displayed;

    \item a \emph{decomposition diagram}, which records the qualitative arrangement of the resonance tongues;

	\item a \emph{boundary diagram}, which plots the rotation number~$\rho$ along~$\mathrm{T}$ and~$\mathrm{T}_0$ as a function of~$\mu$, together with its values at~$\mathrm{H}_1$ and~$\mathrm{H}_2$. The result is a closed curve in the $(\mu,\rho)$-plane that represents~$B$, and each resonance tongue is shown as a horizontal line segment connecting two distinct points on this curve at the corresponding rational value of~$\rho$.
\end{enumerate}

Figure~\ref{bd:fig1}(a) for~$\eta = -0.214$ shows that~$\rho$ varies monotonically along any path from~$\mathrm{H}_1$ to~$\mathrm{H}_2$ that is a graph over the~$\mu$-axis; this is reflected in the continuous color transition from dark blue to dark red as~$\mu$ increases from~$\mathrm{H}_1$ to~$\mathrm{H}_2$. The same monotonicity holds along each of the components~$\mathrm{T}_0$ and~$\mathrm{T}$ of~$B$, so that~$\mathrm{H}_1$ and~$\mathrm{H}_2$ are, respectively, the global minimum and maximum of~$\rho$. Along~$\mathrm{T}$, the rotation number is determined by the Floquet multipliers of the bifurcating periodic solution; along~$\mathrm{T}_0$, it is determined by the frequency of the underlying periodic solution.

It follows that each horizontal line segment in the boundary diagram of panel~(c) has one endpoint on~$\mathrm{T}_0$ and the other on~$\mathrm{T}$, corresponding in panel~(a) to a resonance tongue that connects the two boundary components. All resonance tongues in case~$\mathrm{I}$ therefore belong to a single family with this property, indicated by the uniform shading in panel~(b).

\begin{figure*}[t!]
  \centering
  \includegraphics{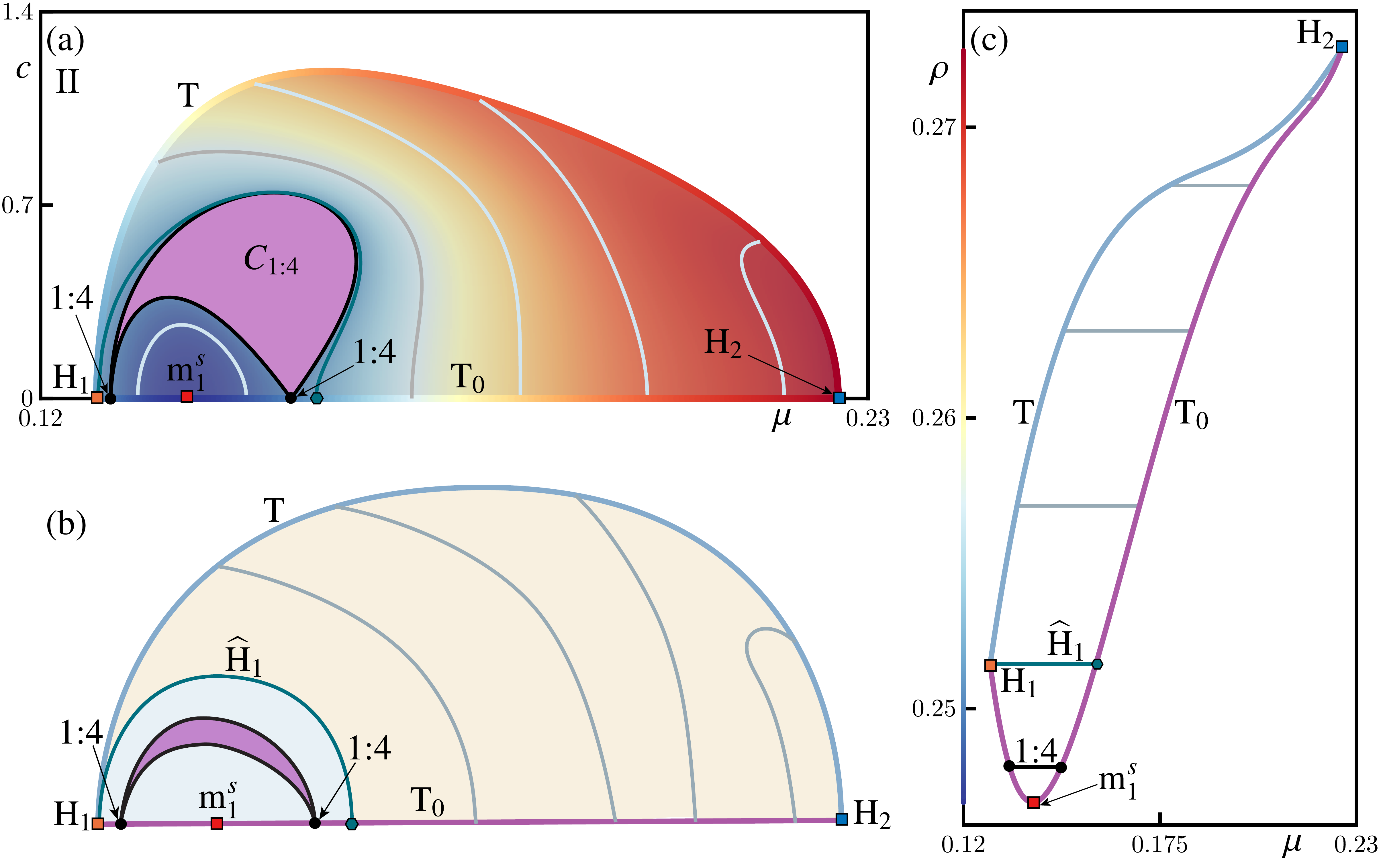}
\caption{\label{bd:fig2} The resonance structure for case~$\mathrm{II}$ at $\eta=-0.236$, in the format of Figure~\ref{bd:fig1}. Shown additionally are the boundary minimum~$\mathrm{m}_1^s$ (red square) and the separatrix~$\widehat{\mathrm{H}}_1$ (turquoise curve). Also shown is a $1{:}4$ resonance tongue, denoted~$C_{1{:}4}$, indicated with black saddle-node boundaries and a purple interior. Panel~(b) shows an additional resonance tongue family in blue, separated by~$\widehat{\mathrm{H}}_1$.}
\end{figure*}
\begin{figure*}[t!]
  \centering
  \includegraphics{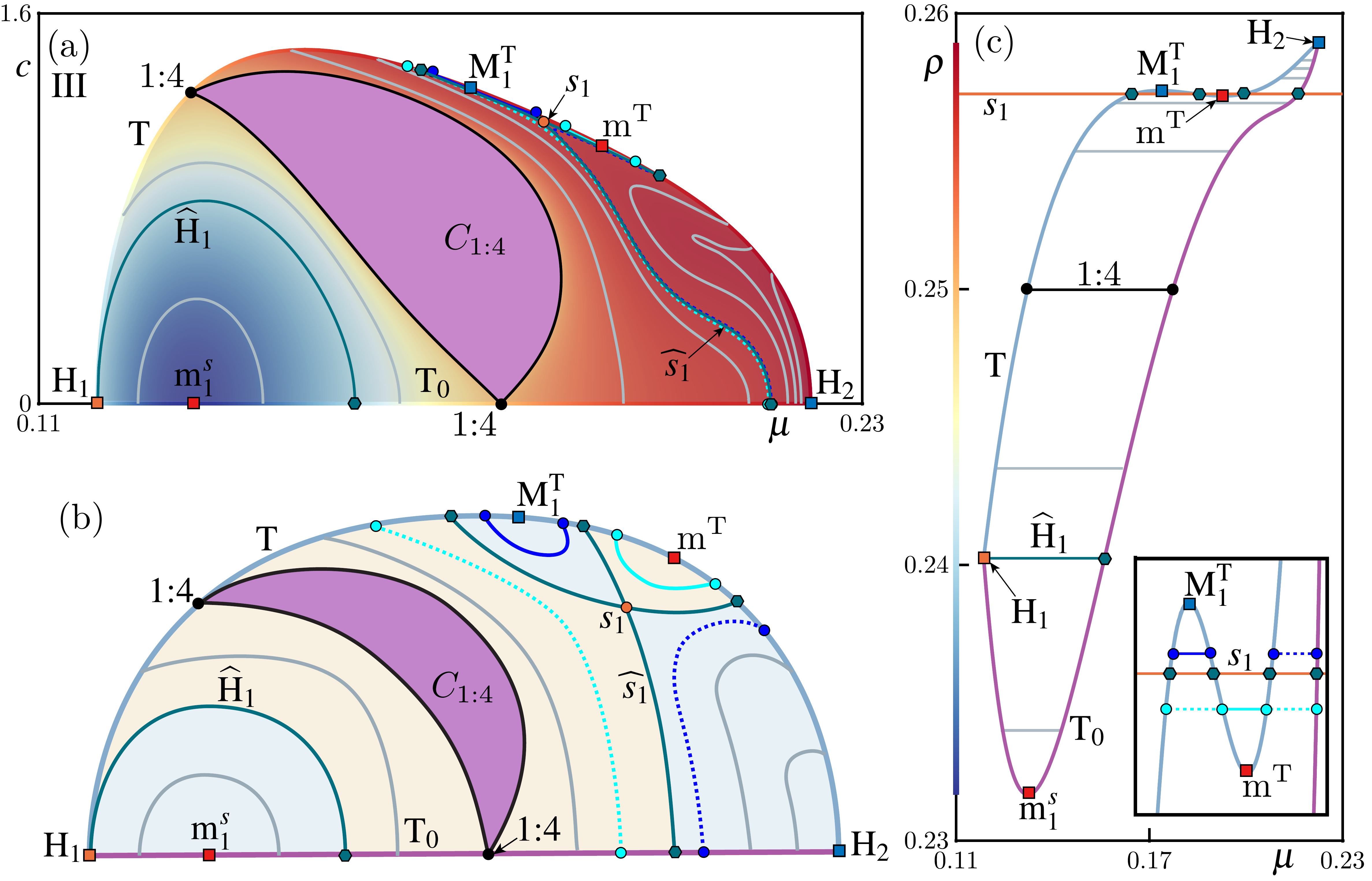}
\caption{\label{bd:fig4} The resonance structure for case~$\mathrm{III}$ at $\eta=-0.2705$, in the format of Figures~\ref{bd:fig1}--\ref{bd:fig2}. Additionally shown are the minimum~$\mathrm{m}^{T}$ (red square) and the maximum~$\mathrm{M}_1^{T}$ (blue square) on the curve~$\mathrm{T}$, together with the saddle~$s_1$ (orange circle) and its separatrix~$\widehat{s}_1$ (turquoise curve) in the interior. Two selected pairs of resonance tongues at distinct values of~$\rho$ near~$s_1$ are shown in cyan and blue, each with a solid and a dotted branch. Panel~(b) now shows the five families of resonance tongues. Panel~(c) contains an inset that enlarges a region near~$\mathrm{m}^{T}$ and~$\mathrm{M}_1^{T}$; the $\rho$-value of~$s_1$ is indicated by the orange line.}
\end{figure*}

In case~$\mathrm{II}$, shown in Figure~\ref{bd:fig2} for~$\eta = -0.236$, the corner point~$\mathrm{H}_1$ is no longer where~$\rho$ attains its minimum: $\rho$ now has a local minimum~$\mathrm{m}_1^s$ on~$\mathrm{T}_0$ and, hence, is no longer monotone on this part of the boundary~$B$. Panel~(c) illustrates this and shows that horizontal lines with values of~$\rho$ near~$\mathrm{m}_1^s$ intersect~$\mathrm{T}_0$ at two points, one on either side of~$\mathrm{m}_1^s$. In Figures~\ref{bd:fig2}(a) and~(b), each such pair of intersections corresponds to a resonance tongue with both endpoints on~$\mathrm{T}_0$; together, these resonance tongues form a new family that includes the prominent $1{:}4$ resonance tongue~$C_{1{:}4}$. Alongside this new family, the family from case~$\mathrm{I}$ persists, with each resonance tongue having one endpoint on~$\mathrm{T}_0$ and the other on~$\mathrm{T}$. The two families are separated in panels~(a) and~(b) by the separatrix~$\widehat{\mathrm{H}}_1$ --- itself either a resonance tongue or a quasi-periodic arc --- which terminates on~$\mathrm{T}_0$ at a point with the same value of~$\rho$ as~$\mathrm{H}_1$; in panel~(c), it corresponds to the horizontal line through~$\mathrm{H}_1$.

In case~$\mathrm{III}$, obtained for the smaller value~$\eta = -0.2705$ and shown in Figure~\ref{bd:fig4}, we find a different resonance tongue configuration. The $1{:}4$ resonance tongue~$C_{1{:}4}$ is again present, but now extends across a larger part of the resonance diagram from a point on~$\mathrm{T}_0$ to a point on~$\mathrm{T}$; however, this does not constitute a topological change of the resonance structure. What distinguishes case~$\mathrm{III}$ is the appearance on the torus bifurcation curve~$\mathrm{T}$ of a maximum~$\mathrm{M}_1^{T}$ and a minimum~$\mathrm{m}^{T}$ of~$\rho$, together with a saddle-type point~$s_1$ in the interior of~$D$. Panels~(a) and~(b) show that~$s_1$ has a separatrix~$\widehat{s}_1$ that encloses narrow regions around~$\mathrm{M}_1^{T}$ and~$\mathrm{m}^{T}$, thus separating these extrema from the rest of~$D$. Since~$s_1$ lies neither on~$\mathrm{T}$ nor on~$\mathrm{T}_0$, it does not appear in the boundary diagram in panel~(c); we indicate its $\rho$-value by a horizontal line. The solid-blue and solid-cyan horizontal segments in panel~(c) represent the resonance tongues that form around~$\mathrm{M}_1^{T}$ and~$\mathrm{m}^{T}$ in panels~(a) and~(b), each with two endpoints on~$\mathrm{T}$, one on either side of the corresponding extremum. The dotted-cyan and dotted-blue segments are at the same $\rho$-values as their solid counterparts, and they represent resonance tongues on the other side of the separatrix~$\widehat{s}_1$, each connecting~$\mathrm{T}_0$ to~$\mathrm{T}$. Overall, there are now five distinct families of resonance tongues, distinguished by the way their endpoints attach to the boundary components~$\mathrm{T}_0$ and~$\mathrm{T}$; these families are separated in the half-disk~$D$ by the separatrices~$\widehat{\mathrm{H}}_1$ and~$\widehat{s}_1$.

The change in the resonance structure from one case to the next is the result of a change on the boundary~$B$: the appearance of the minimum~$\mathrm{m}_1^s$ on~$\mathrm{T}_0$ for the transition to case~$\mathrm{II}$, and of the maximum~$\mathrm{M}_1^{T}$ and the minimum~$\mathrm{m}^{T}$ on~$\mathrm{T}$ for the transition to case~$\mathrm{III}$. Specifically, in case~$\mathrm{II}$, the corner point~$\mathrm{H}_1$ acquires the separatrix~$\widehat{\mathrm{H}}_1$, which partitions the resonance tongues into two families. In case~$\mathrm{III}$, the change on the boundary is accompanied by the appearance of an interior saddle~$s_1$, whose separatrix~$\widehat{s}_1$ further partitions the resonance tongues. Cases~$\mathrm{I}$--$\mathrm{III}$ therefore suggest that the organization of resonance tongues in~$D$ is governed by the singularities of~$\rho$, both on the boundary~$B$ and in the interior of~$D$. In the next section, we formalize this viewpoint and develop the tools needed to classify such configurations and the transitions between them.

\section{Topological framework}
\label{section:topo_framework}
We now present the topological framework for describing the global arrangement of resonance tongues. For fixed~$\eta$, the rotation number is a continuous scalar field
\[
    \rho_\eta \colon D_\eta \to [0,1]
\]
on the half-disk region~$D_\eta$ and extends continuously to the boundary~$B_\eta = \mathrm{T} \cup \mathrm{T}_0 \cup \{\mathrm{H}_1, \mathrm{H}_2\}$, where it is smooth over the curves~$\mathrm{T}$ and~$\mathrm{T}_0$. The \emph{resonance surface} is the graph of~$\rho_\eta$ over~$D_\eta$,
\begin{equation}
  \mathcal{S}_\eta
  =
  \{(\mu,c,\rho_\eta(\mu,c)) \mid (\mu,c)\in D_\eta\}
  \subset \mathbb{R}^2 \times [0,1],
  \label{eq:resonance_surface}
\end{equation}
whose boundary is the graph of~$\rho_\eta$ over~$B_\eta$. In this representation, resonance tongues appear as level sets of~$\mathcal{S}_\eta$ at rational values of~$\rho_\eta$, and their global arrangement is encoded by the level-set topology of~$\mathcal{S}_\eta$. Unless stated otherwise, we use the same notation for an object on~$\mathcal{S}_\eta$ and for its projection to~$D_\eta$.

As any periodically forced system with sinusoidal forcing, system~\eqref{eq:forced_ODE_second} is invariant under the reflection~$c \mapsto -c$ combined with the half-period phase shift~$\theta \mapsto \theta + \tfrac{1}{2}$; this conjugacy preserves the rotation number, so the symmetry carries over to the resonance surface~$\mathcal{S}_\eta$, extended by reflection to~$c < 0$, as the involution
\begin{equation}
    (\mu,c,\rho) \mapsto (\mu,-c,\rho),
    \label{eq:involution}
\end{equation}
whose fixed-point set is the zero-forcing segment~$\mathrm{T}_0$, which we also refer to as the \emph{symmetry segment}.

\begin{figure}[t!]
  \centering
  \includegraphics{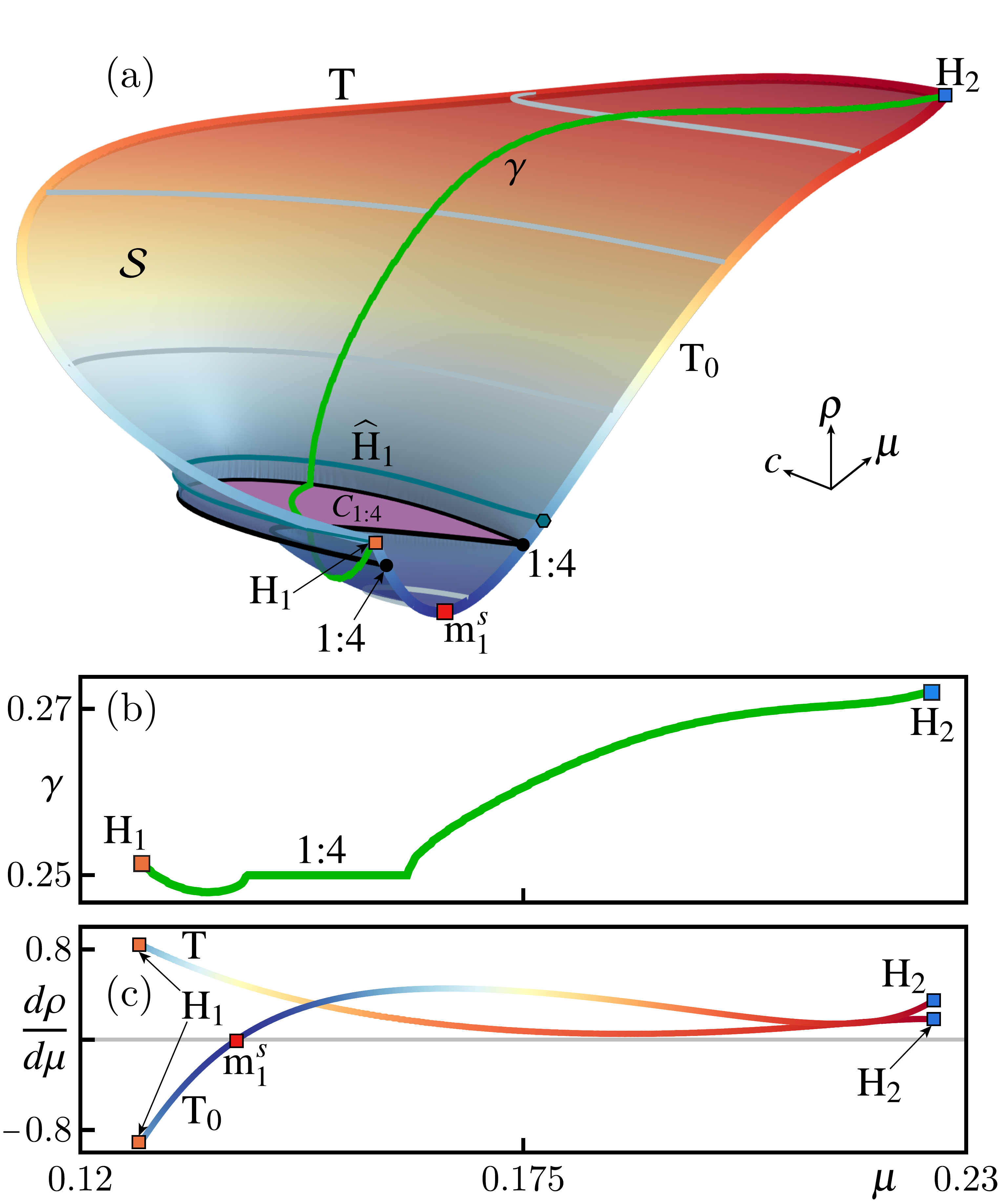}
\caption{\label{fig:intro_theory_fig} Resonance surface~$\mathcal{S}$ of case~$\mathrm{II}$ with~$\eta = -0.236$, showing the same objects as in Figure~\ref{bd:fig2}. Panel~(a) shows~$\mathcal{S}$ in $(\mu,c,\rho)$-space, with the path~$\gamma$ (green curve) from~$\mathrm{H}_1$ to~$\mathrm{H}_2$. Panel~(b) shows the rotation number~$\rho$ along~$\gamma$ as a function of~$\mu$, and panel~(c) shows the derivative~$\frac{d\rho}{d\mu}$ along the boundary curves~$\mathrm{T}$ and~$\mathrm{T}_0$ as a function of~$\mu$.}
\end{figure}

We revisit case~$\mathrm{II}$ from Figure~\ref{bd:fig2} to illustrate why the standard Morse-theoretic notion of an isolated critical point does not suffice for describing~$\mathcal{S}$: there is a lack of smoothness, both in the interior of~$\mathcal{S}$ and on its boundary~$B$. The rendering of~$\mathcal{S}$ as a surface in~$(\mu,c,\rho)$-space in Figure~\ref{fig:intro_theory_fig}(a) shows the~$1{:}4$ resonance tongue as a sizable region $C_{1{:}4}$ of nonzero area on which~$\rho$ is constant. Panel~(b) shows~$\rho$ along the smooth path~$\gamma$ from~$\mathrm{H}_1$ to~$\mathrm{H}_2$, with a corresponding constant portion that crosses~$C_{1{:}4}$. Indeed, this one-dimensional cross section of~$\mathcal{S}$ is a devil's staircase, which is constant at every rational value of~$\rho$; however, all other `steps' are too small to distinguish in Figure~\ref{fig:intro_theory_fig}(b). On the segments~$\mathrm{T}$ and~$\mathrm{T}_0$ of the boundary~$B$, the rotation number~$\rho$ is smooth, but this is not the case at the corner points~$\mathrm{H}_1$ and~$\mathrm{H}_2$. This is illustrated in panel~(c), where the derivative~$\frac{d\rho}{d\mu}$ along~$\mathrm{T}$ and~$\mathrm{T}_0$ has a jump discontinuity at~$\mathrm{H}_1$ and~$\mathrm{H}_2$. We conclude that a useful notion of singularity for~$\mathcal{S}$ must accommodate regions of nonzero area in the interior, critical points on the smooth parts~$\mathrm{T}$ and~$\mathrm{T}_0$ of~$B$, and the corner points~$\mathrm{H}_1$ and~$\mathrm{H}_2$. We now proceed to formalize this observation in purely topological terms.

\subsection{Regular and singular contours of~$\mathcal{S}$}
\label{subsec:sing_res_surf}
We now introduce, for fixed~$\eta$, a single notion of singularity that encompasses all parts of a resonance surface~$\mathcal{S}$ over the compact half-disk~$D$ with boundary~$B$. The idea is to track changes in the connected components of the sub- and super-level sets of~$\mathcal{S}$ across individual contours.

The level set of~$\mathcal{S}$ at a given height~$h$ is
\begin{equation}
  L_{h} \;=\;
  \{\,(\mu,c,\rho(\mu,c))\in\mathcal{S} \mid \rho(\mu,c)=h\,\},
\end{equation}
and a connected component of~$L_h$ is called a \emph{contour}. To monitor how the level sets are arranged as~$h$ sweeps across~$[0,1]$, we collect them into two one-parameter collections. The \emph{sub-level set}
\begin{equation}
  \mathcal{L}^-_{h}
  \;=\;
  \{\, (\mu,c,\rho(\mu,c)) \in \mathcal{S} \;\mid\; \rho(\mu,c)\le h \,\}
\end{equation}
contains all points of~$\mathcal{S}$ with rotation number at most~$h$, and the \emph{super-level set}
\begin{equation}
  \mathcal{L}^+_{h}
  \;=\;
  \{\, (\mu,c,\rho(\mu,c)) \in \mathcal{S} \;\mid\; \rho(\mu,c)\ge h \,\}
\end{equation}
contains those with rotation number at least~$h$. As~$h$ sweeps across~$[0,1]$, the collection~$\{\mathcal{L}^-_h\}$ grows monotonically from~$h = \min\{\rho\}$ to~$h = \max\{\rho\}$, while~$\{\mathcal{L}^+_h\}$ shrinks monotonically; each forms a continuous filtration of~$\mathcal{S}$~\cite{edelsbrunner2010computational, matsumoto2002introduction}.

\begin{figure*}[t!]
  \centering
  \includegraphics[width=\linewidth]{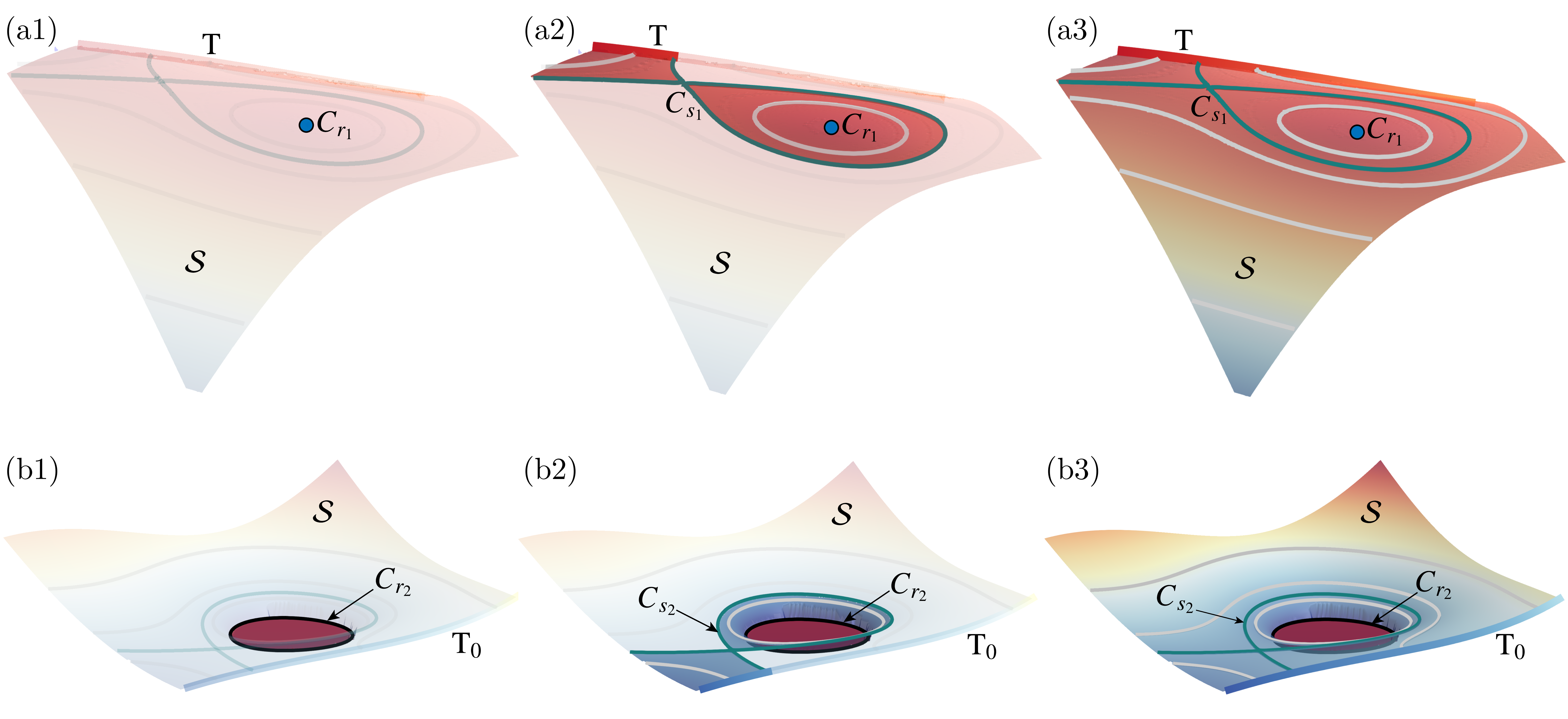}
\caption{Sub- and super-level sets of the resonance surface~$\mathcal{S}$. Panels~(a1)--(a3), for $\eta = -0.286$, show the super-level sets~$\mathcal{L}^{+}_{h}$ for three decreasing values of~$h$, with the interior maximum~$C_{r_1}$ (blue dot) and the monotone saddle~$C_{s_1}$ (turquoise contour). Panels~(b1)--(b3), for $\eta = -0.277$, show the sub-level sets~$\mathcal{L}^{-}_{h}$ for increasing~$h$, with the monotone saddle~$C_{s_2}$ and the interior minimum~$C_{r_2}$.}
\label{fig:example_filtration}
\end{figure*}

The number of connected components of either sub- or super-level set is a topological invariant, the zeroth Betti number, denoted~$\beta_0(\mathcal{L}^-_h)$ and~$\beta_0(\mathcal{L}^+_h)$ respectively~\cite{lee2000introduction}. To detect when these counts change and to attribute the change to an individual contour~$C$, we introduce the up- and down-degrees of~$C$, adapting analogous notions from computational topology~\cite{edelsbrunner2010computational, carr2003computing}.
\begin{definition}[Up- and down-degree of a contour]
  \label{def:degree_value}
  Fix~$h \in [0,1]$ and let~$C \subset L_h$ be a contour of~$\mathcal{S}$. For~$\delta > 0$, let
  \begin{equation*}
      U_\delta \;=\; \bigl\{\, p \in \mathcal{S} \;\big|\; \|p - q\| < \delta
      \text{ for some } q \in C \,\bigr\}
  \end{equation*}
  be the open~$\delta$-neighborhood of~$C$ in~$\mathcal{S}$ with respect to the Euclidean metric on~$\mathbb{R}^3$, and define the local component counts
  \begin{equation*}
    \beta^{-}_\delta(h) \;=\; \beta_0\!\bigl(\mathcal{L}^-_h \cap U_\delta\bigr),
    \qquad
    \beta^{+}_\delta(h) \;=\; \beta_0\!\bigl(\mathcal{L}^+_h \cap U_\delta\bigr).
  \end{equation*}
  Assume there exist~$\delta_0 > 0$ and~$\varepsilon_0 > 0$ such that~$\beta^\pm_\delta(h \pm \varepsilon)$ is constant for all~$0 < \delta < \delta_0$ and~$0 < \varepsilon < \varepsilon_0$. The \emph{down-degree} of~$C$ is
  \begin{equation*}
    \deg^{-}(C) \;=\; \beta^{-}_\delta(h - \varepsilon),
  \end{equation*}
  and the \emph{up-degree} of~$C$ is
  \begin{equation*}
    \deg^{+}(C) \;=\; \beta^{+}_\delta(h + \varepsilon).
  \end{equation*}
\end{definition}

Each degree is a non-negative integer counting the connected components of the corresponding filtration immediately below or above~$C$. A regular contour has one such component on each side, while higher counts indicate the birth or merging of components across~$C$.

\begin{definition}[Regular and singular contours of~$\mathcal{S}$]
  \label{def:regular_singular_contour}
  Fix a height~$h \in [0,1]$ such that the level set~$L_h$ contains no corner points, that is,~$\mathrm{H}_1 \notin L_h$ and~$\mathrm{H}_2 \notin L_h$. Let~$C \subset L_h$ be a contour. We call~$C$ \emph{regular} if
  \[
    \deg^{-}(C) = \deg^{+}(C) = 1,
  \]
  and \emph{singular} otherwise. Contours that contain a corner point are regarded as singular by convention.
\end{definition}

The corner-point convention is needed because such contours may organize the arrangement of nearby resonance tongues without producing a change in the number of connected components of either filtration. For example, in Figure~\ref{fig:intro_theory_fig}, the corner maximum~$\mathrm{H}_2$ has up-degree~$0$ and is therefore already detected as singular by its degree count (see Definition~\ref{def:regular_singular_contour}), whereas the corner saddle~$\mathrm{H}_1$ has up- and down-degrees both equal to one --- so the filtrations register no change there --- yet it still plays an important organizing role. Our convention ensures that both are classified as singular.

We now distinguish singular contours by their location on~$\mathcal{S}$, separating those that lie in the interior from those that meet the boundary~$B$ either at a critical point of~$\rho$ or at a corner point.
\begin{definition}[Corner, boundary, and interior singularities]
\label{def:boundary_vs_interior}
For $h\in[0,1]$, let $C \subset L_{h}$ be a singular contour of~$\mathcal{S}$. We call~$C$:
\begin{enumerate}[label=(\roman*)]
\item a \emph{corner singularity} if it contains a corner point;
\item a \emph{boundary singularity} if it contains a critical point of~$\rho$ on either the torus bifurcation curve~$\mathrm{T}$ or the symmetry segment~$\mathrm{T}_0$;
\item an \emph{interior singularity} if neither (i) nor (ii) holds.
\end{enumerate}
\end{definition}

Boundary and interior singularities are further classified as maximum, minimum, and saddle-type by their down- and up-degrees; corner points are classified separately.
\begin{definition}[Classification of boundary and interior singularities]
\label{def:local_singularity_type}
Let~$C \subset L_h$ be an boundary or interior singular contour with down- and up-degrees~$\deg^{-}(C)$ and~$\deg^{+}(C)$. We distinguish the following types:
\begin{enumerate}[label=(\roman*)]
  \item \emph{minimum} when~$\deg^{-}(C) = 0$ and $\deg^{+}(C) = 1$;
  \item \emph{maximum} when~$\deg^{-}(C) = 1$ and $\deg^{+}(C) = 0$;
  \item \emph{monotone saddle} when~$\deg^{-}(C) = 2$ and $\deg^{+}(C) = 1$, or $\deg^{-}(C) = 1$ and $\deg^{+}(C) = 2$;
  \item \emph{compound saddle of order}~$(m,k)$ \emph{with}~$m,k \geq 2$ when~$\deg^{-}(C) = m$ and $\deg^{+}(C) = k$.
\end{enumerate}
\end{definition}
Interior minima, interior maxima, and monotone saddles are direct analogues of the generic critical points of classical Morse theory, with Morse indices~$0$,~$2$, and~$1$ respectively~\cite{audin2014morse}. A compound saddle is non-generic in the smooth setting, where it corresponds to two or more saddle-type critical points sharing a common level set. In our setting, on the other hand, compound saddles are generic: the piecewise-constant structure of~$\rho$ on~$\mathcal{S}$, as well as its boundary, allows several such `merge events' to occur simultaneously along a single contour.

Figure~\ref{fig:example_filtration} illustrates how the filtrations~$\{\mathcal{L}^-_h\}$ and~$\{\mathcal{L}^+_h\}$ identify interior singularities. Panels~(a1)--(a3) display the super-level set~$\mathcal{L}^{+}_{h}$ at three decreasing values of~$h$. As~$h$ decreases past the level of~$C_{r_1}$, a new component appears, so~$C_{r_1}$ is an interior maximum; as~$h$ decreases further past~$C_{s_1}$, two components merge, so~$C_{s_1}$ is a monotone saddle. Panel~(a3) shows the single remaining component. Panels~(b1)--(b3) display the analogous picture for the sub-level set~$\mathcal{L}^{-}_{h}$ at three increasing values of~$h$. As~$h$ increases past~$C_{r_2}$, a new component appears, so~$C_{r_2}$ is an interior minimum; as~$h$ increases further past~$C_{s_2}$, two components merge, so~$C_{s_2}$ is a monotone saddle. Panel~(b3) shows the single remaining component. The difference between panels~(a1)--(a3) and panels~(b1)--(b3) is that~$C_{r_2}$ is clearly a large plateau while~$C_{r_1}$,~$C_{s_1}$, and~$C_{s_2}$ are very thin; the filtration detects all four in exactly the same way.

As we discussed, corner points are special and are classified according to the monotonicity of~$\rho$ along the curves~$\mathrm{T}$ and~$\mathrm{T}_0$.
\begin{definition}[Corner extrema and corner saddles]
\label{def:corner_classification}
Let~$\mathrm{H}$ be one of the corner points~$\mathrm{H}_1$ or~$\mathrm{H}_2$. We classify~$\mathrm{H}$ according to the monotonicity of~$\rho$ along~$\mathrm{T}$ and~$\mathrm{T}_0$ as one moves away from~$\mathrm{H}$:
\begin{enumerate}[label=(\roman*)]
    \item $\mathrm{H}$ is a \emph{corner minimum} if~$\rho$ is increasing along both~$\mathrm{T}$ and~$\mathrm{T}_0$;
    \item $\mathrm{H}$ is a \emph{corner maximum} if~$\rho$ is decreasing along both~$\mathrm{T}$ and~$\mathrm{T}_0$;
    \item $\mathrm{H}$ is a \emph{corner saddle} otherwise.
\end{enumerate}
\end{definition}

\subsection{Geometry of singularities}
\label{subsec:boundary}
\begin{figure}[ht!]
  \centering
  \includegraphics[width=0.5\linewidth]{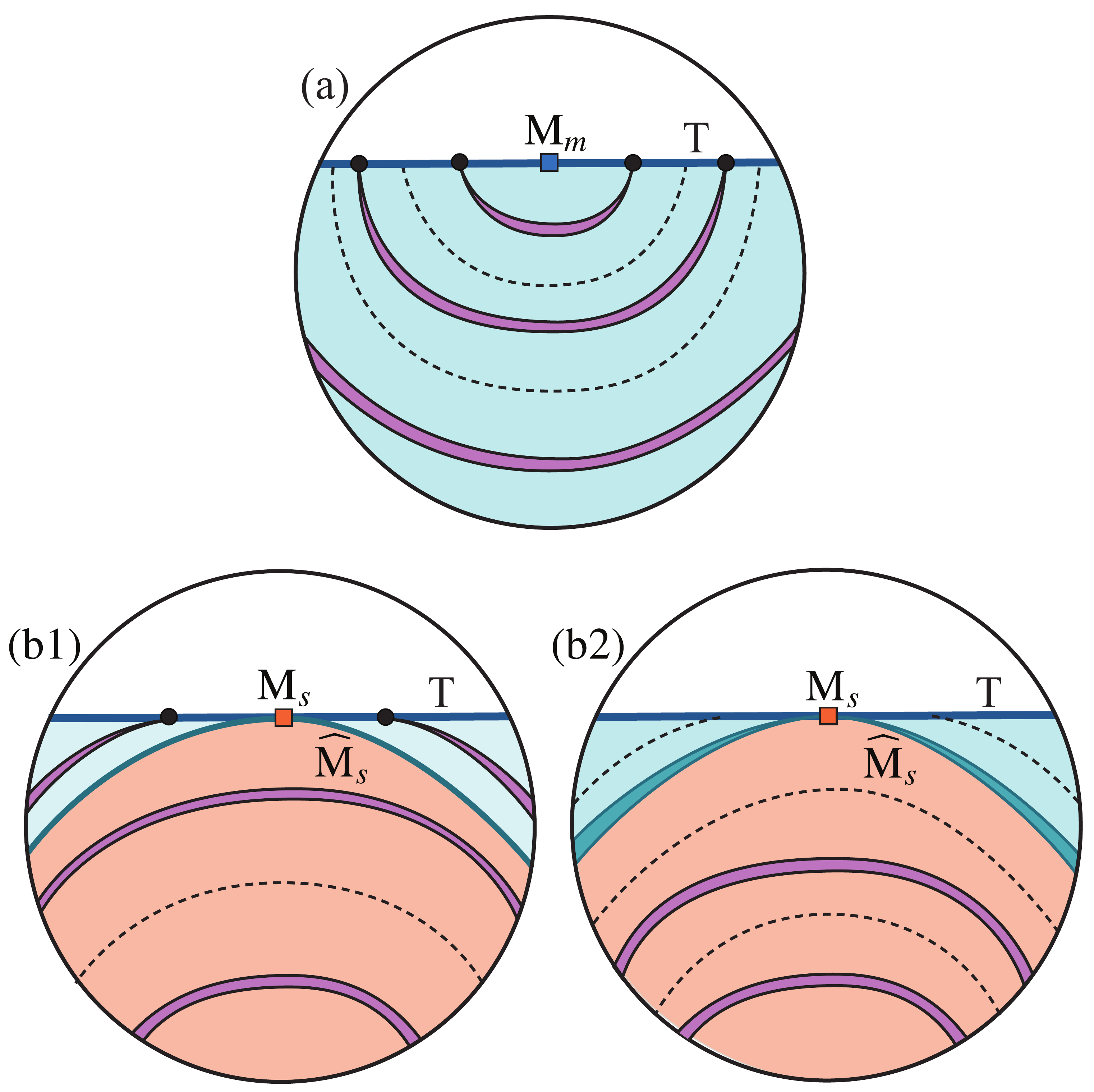}
\caption{Local decomposition diagrams of boundary singularities on~$\mathrm{T}$. Blue shading indicates regions of~$\mathcal{S}$ below the singularity, and orange shading regions above it. Resonance tongues are shown in purple with black boundaries, and quasi-periodic arcs as dotted curves. Panel~(a) shows a boundary maximum~$\mathrm{M}_{m}$ (blue square); panels~(b1) and~(b2) show a boundary saddle~$\mathrm{M}_{s}$ (orange square) with its separatrix~$\widehat{\mathrm{M}}_{s}$ (dark turquoise) for irrational and rational~$\rho$, respectively.}
\label{fig:topo_sketches_T}
\end{figure}

We have so far defined and classified the singularities of the resonance surface~$\mathcal{S}$ by their degrees and their location. We now describe the local arrangement of regular contours --- equivalently, of resonance tongues and quasi-periodic arcs --- around each singularity. Saddles and interior extrema, in particular, are treated separately for rational and irrational rotation number~$\rho$. We present each singularity through a decomposition diagram of its neighborhood in the~$(\mu,c)$-plane, where we distinguish between families of resonance tongues that lie above and below the associated singular contour. For interior singularities, we additionally display the surface~$\mathcal{S}$ itself. Throughout, we present one representative of each singularity type; the remaining configurations are obtained by replacing~$\rho$ with~$-\rho$.

\subsubsection{Boundary singularities}
The arrangement of families of resonance tongues around a boundary singularity depends on whether the singularity lies on the torus bifurcation curve~$\mathrm{T}$ or on the symmetry segment~$\mathrm{T}_0$. Along~$\mathrm{T}$, the resonance surface~$\mathcal{S}$ ends and so the decomposition diagram resembles a half-neighborhood of the singularity. Along~$\mathrm{T}_0$, on the other hand, the decomposition diagram is invariant under the reflection~$c \mapsto -c$, and~$\mathrm{T}_0$ acts as a symmetry axis joining two reflected halves; for clarity, we display all relevant objects for both signs of~$c$ at singularities on~$\mathrm{T}_0$ and at corner points. For definiteness, we consider a local maximum~$\mathrm{M}_m$ of the rotation number~$\rho$ when~$\rho$ is restricted to the boundary~$B$; the case of a boundary minimum of~$\rho$ on~$B$ follows by considering~$-\rho$. 

\begin{figure}[t!]
  \centering
  \includegraphics[width=0.5\linewidth]{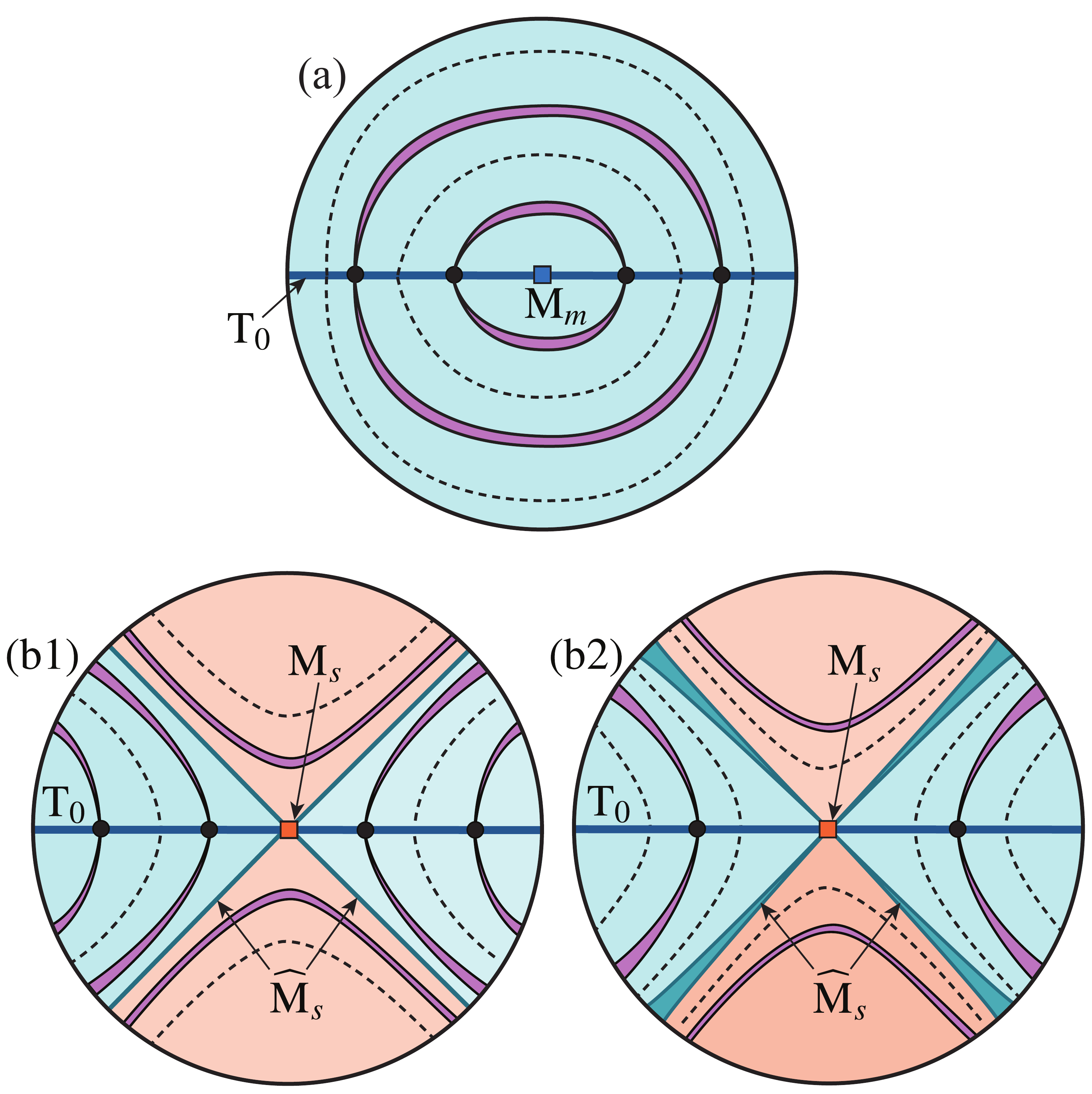}
\caption{Local decomposition diagrams of boundary singularities on~$\mathrm{T}_0$, presented as in Figure~\ref{fig:topo_sketches_T}. Panel~(a) shows a boundary maximum~$\mathrm{M}_{m}$; panels~(b1) and~(b2) show a boundary saddle~$\mathrm{M}_{s}$ with its separatrix~$\widehat{\mathrm{M}}_{s}$ for irrational and rational~$\rho$, respectively.}
\label{fig:topo_sketches_T0}
\end{figure}

Figure~\ref{fig:topo_sketches_T} illustrates the local decompositions for a boundary maximum and a boundary saddle on~$\mathrm{T}$. Panel~(a) shows that near a maximum~$\mathrm{M}_m$, there is a single family of regular contours, each with an endpoint on~$\mathrm{T}$ on either side of~$\mathrm{M}_m$. Row~(b) shows a boundary saddle~$\mathrm{M}_s$ together with its separatrix~$\widehat{\mathrm{M}}_s$, which comprises two branches that meet tangentially at~$\mathrm{M}_s$. When~$\rho$ is irrational, as in panel~(b1), $\widehat{\mathrm{M}}_s$ consists of two arcs; when~$\rho$ is rational, as in panel~(b2), $\widehat{\mathrm{M}}_s$ consists of a pair of resonance tongues with both endpoints at~$\mathrm{M}_s$. Regardless of the rationality of~$\rho$, the singular contour~$\widehat{\mathrm{M}}_s$ partitions nearby regular contours into three families: two in the blue regions, each meeting~$\mathrm{T}$ at a single point, and a third in the red region with no endpoints near~$\mathrm{M}_s$. These configurations have been observed previously in the context of varying a third parameter~\cite{terrien2023merging}, where resonance tongues were found to emerge at a boundary maximum and to disconnect at a boundary saddle on~$\mathrm{T}$; however, the rational and irrational sub-cases were not distinguished there.

\begin{figure}[ht!]
  \centering
  \includegraphics[width=0.5\linewidth]{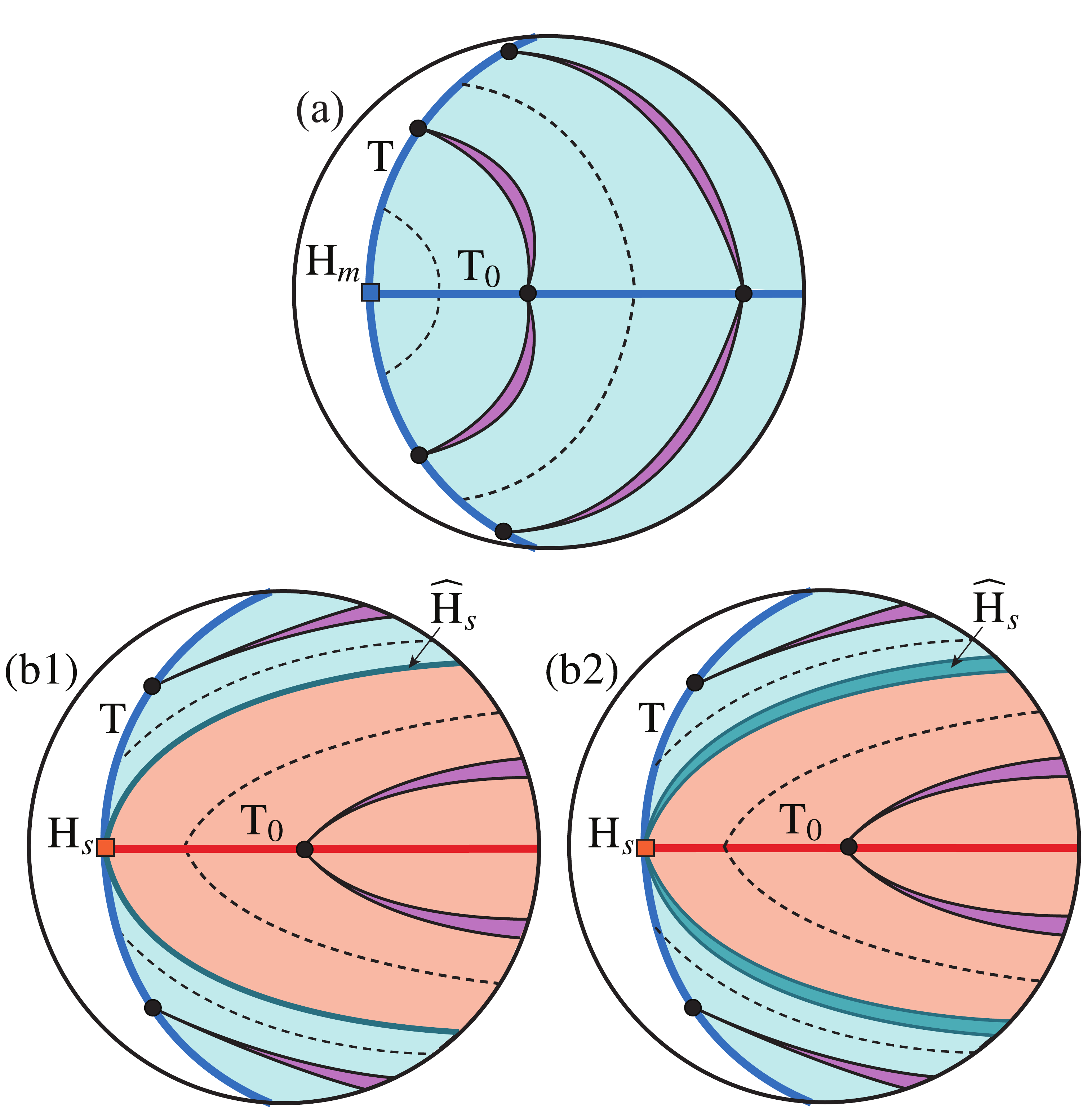}
\caption{Local decomposition diagrams of corner singularities presented as in Figures~\ref{fig:topo_sketches_T} and~\ref{fig:topo_sketches_T0}. Panel~(a) shows a corner maximum~$\mathrm{H}_{m}$; panels~(b1) and~(b2) show a corner saddle~$\mathrm{H}_{s}$ with separatrix~$\widehat{\mathrm{H}}_{s}$ for irrational and rational number~$\rho$, respectively.}
\label{fig:topo_sketches_H}
\end{figure}

Figure~\ref{fig:topo_sketches_T0} shows the decomposition around a boundary maximum and a boundary saddle on~$\mathrm{T}_0$; the arrangement of resonance tongues resembles that in Figure~\ref{fig:topo_sketches_T}, but is extended to both signs of~$c$ by reflection. Panel~(a) shows a single family of regular contours encircling the boundary maximum~$\mathrm{M}_m$, each ending on~$\mathrm{T}_0$ at two points, one on either side of~$\mathrm{M}_m$. Row~(b) shows a boundary saddle~$\mathrm{M}_s$ together with its separatrix~$\widehat{\mathrm{M}}_s$, which has two branches on each side of~$\mathrm{T}_0$; unlike in the case of a boundary saddle on~$\mathrm{T}$, these branches meet transversally at~$\mathrm{M}_s$. The separatrix~$\widehat{\mathrm{M}}_s$ again partitions the nearby regular contours with~$c>0$ into three families: two in the blue regions, each ending on~$\mathrm{T}_0$ at a single point, and a third in the red region.

Figure~\ref{fig:topo_sketches_H} shows the decomposition near a corner maximum and a corner saddle at the intersection of~$\mathrm{T}_0$ and~$\mathrm{T}$. The local resonance structure resembles that near singularities on~$\mathrm{T}$, with the additional reflection symmetry about~$\mathrm{T}_0$. Panel~(a) shows a single family of regular contours encircling the corner maximum~$\mathrm{H}_m$, each with one endpoint on~$\mathrm{T}_0$ and one endpoint on~$\mathrm{T}$. Row~(b) shows a corner saddle~$\mathrm{H}_s$ together with its separatrix~$\widehat{\mathrm{H}}_s$, which has branches on each side of~$\mathrm{T}_0$ that meet tangentially at~$\mathrm{H}_s$. For both irrational~$\rho$ in panel~(b1) and rational~$\rho$ in panel~(b2), $\widehat{\mathrm{H}}_s$ partitions the nearby regular contours with~$c>0$ into two families: one in the blue region, with an endpoint on~$\mathrm{T}$, and a second in the red region, with an endpoint on~$\mathrm{T}_0$.

\subsubsection{Interior singularities}
\begin{figure}[t!]
  \centering
  \includegraphics{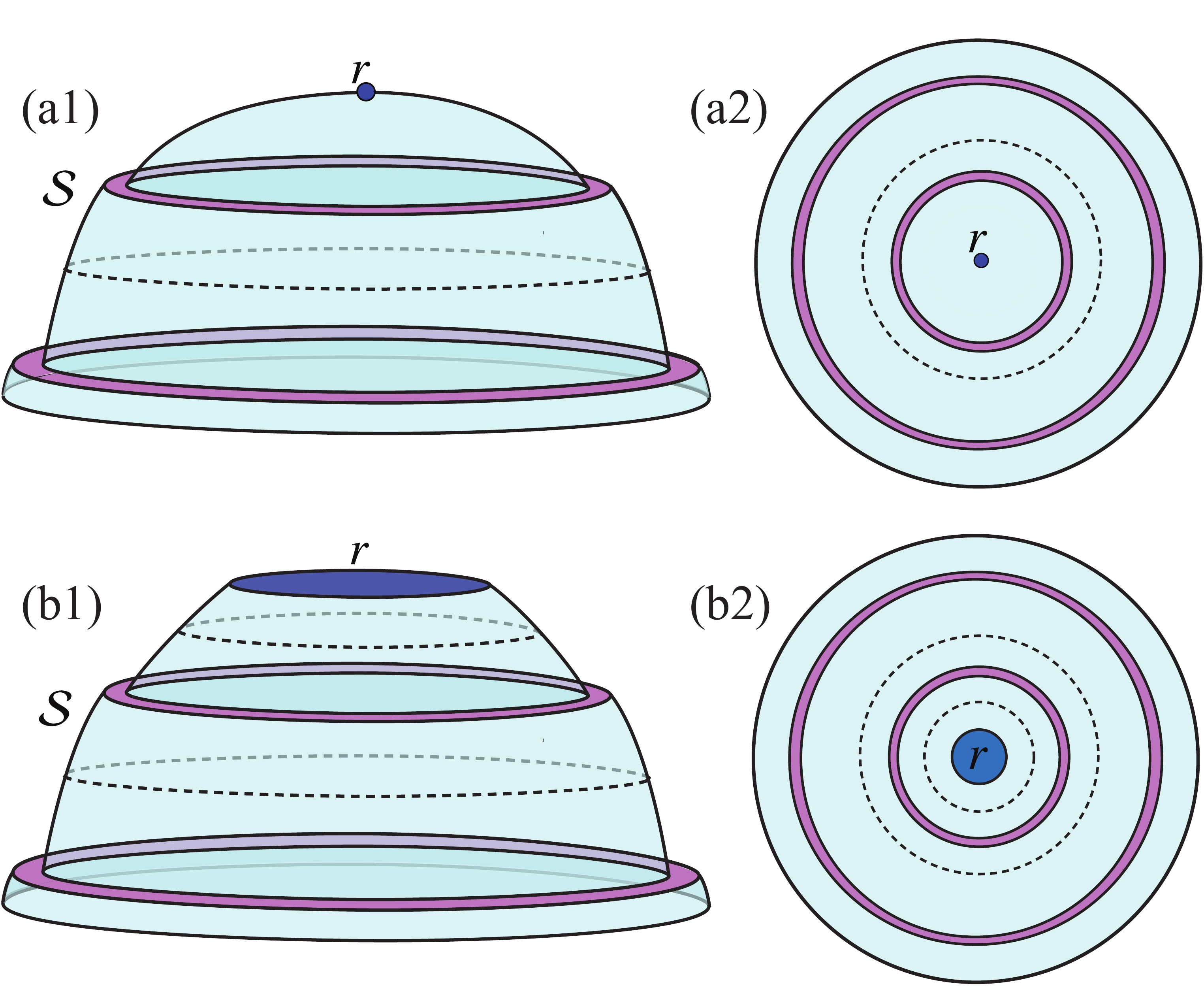}
\caption{Local resonance surfaces and decomposition diagrams near an interior maximum~$r$, with the same color convention as Figure~\ref{fig:topo_sketches_H}. Panel~(a1) shows the resonance surface~$\mathcal{S}$ for irrational~$\rho$ at~$r$ (blue dot) and panel~(a2) its decomposition diagram. Panels~(b1)--(b2) show the same for rational~$\rho$ at~$r$ (blue terrace).}
\label{fig:topo_sketch_3d_max}
\end{figure}

As for a smooth surface, an interior singularity on the resonance surface~$\mathcal{S}$ is generically either a local extremum or a saddle. Figure~\ref{fig:topo_sketch_3d_max} is a sketch of~$\mathcal{S}$ and its decomposition near an interior maximum~$r$. Panels~(a1) and~(b1) show that~$\mathcal{S}$ resembles an inverted bowl lying entirely below~$r$, so its regular contours are closed; in panels~(a2) and~(b2), they appear respectively as one-dimensional loops and annular resonance tongues encircling~$r$. For irrational~$\rho$, as in row~(a), the maximum~$r$ is a point, while for rational~$\rho$, as in row~(b), it is instead a two-dimensional terrace in the form of a plateau. 
  
\begin{figure}[t!]
  \centering  
  \includegraphics{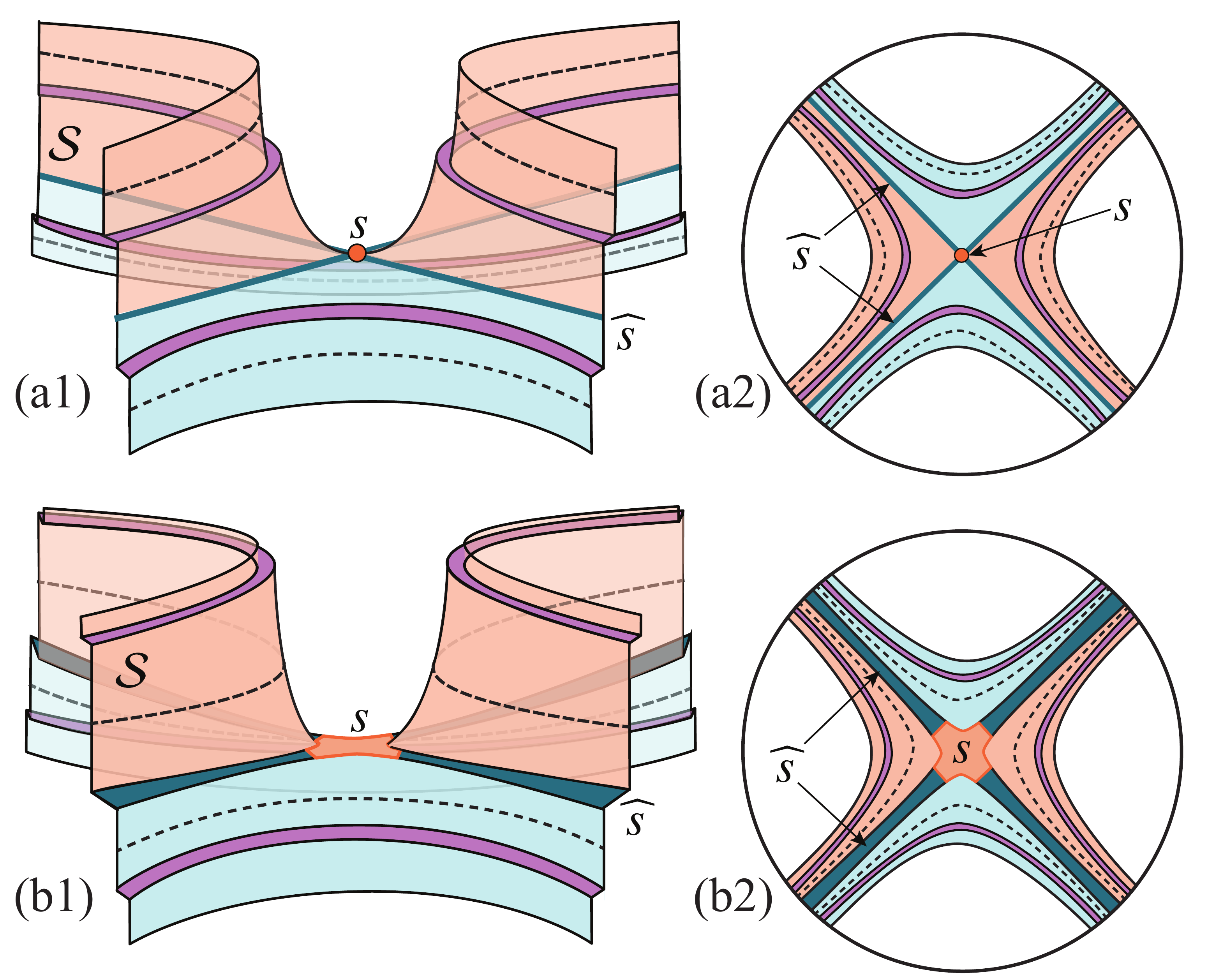}
\caption{Local resonance surfaces and decomposition diagrams near a monotone saddle~$C_{s_1}$, shown in the same manner as Figure~\ref{fig:topo_sketch_3d_max}. Shown is~$C_{s_1}$, decomposed into the inner saddle~$s_1$ (orange) and its separatrices~$\widehat{s}_1$ (turquoise) for irrational~$\rho$ in panels~(a1)--(a2) and rational~$\rho$ in panels~(b1)--(b2).}
\label{fig:topo_sketches3}
\end{figure}

Figure~\ref{fig:topo_sketches3} shows similarly the case of a saddle-type singularity in the form of a plateau~$C_{s_1}$, together with the corresponding local decomposition diagram. This situation is found in system~\eqref{eq:forced_ODE_second} for all monotone saddles and for several~$(2,2)$-compound saddles: the contour~$C_{s_1}$ comprises four separatrix branches, each labeled~$\widehat{s}_1$, that meet at a central object that we refer to colloquially as the \emph{saddle}, denoted~$s_1$. For irrational~$\rho$, as in row~(a), the branches~$\widehat{s}_1$ are one-dimensional arcs that meet at the saddle point~$s_1$; for rational~$\rho$, as in row~(b), both~$s_1$ and~$\widehat{s}_1$ are `fattened' parts of a resonance tongue: the branches~$\widehat{s}_1$ are rectangular strips attached to an inner saddle region~$s_1$; together they form a saddle separatrix contour with four local separatrices. In panels~(a1) and~(b1), the saddle contour $C_{s_1}$ divides~$\mathcal{S}$ into two components above and two below, each containing a separate family of regular contours. The contours of these four families resemble hyperbolae in panels~(a2) and~(b2).

\subsection{Decompositions of the resonance surface}
\label{subsec:decomp_res_surf}
We now introduce a decomposition of the regular contours of the resonance surface~$\mathcal{S}$, by adapting the standard contour equivalence from level-set topology~\cite{edelsbrunner2010computational}. The decomposition groups regular contours by whether they can be connected without crossing a singular contour.
\begin{definition}[Equivalence of regular contours]
\label{def:contour_equivalence}
Let~$\mathcal{C}_{\mathrm{sing}}$ denote the set of all singular contours on~$\mathcal{S}$, and let~$\mathcal{C}_{\mathrm{reg}}$ denote the set of regular contours, as in Definition~\ref{def:regular_singular_contour}. For~$C, C' \in \mathcal{C}_{\mathrm{reg}}$, we say that~$C$ and~$C'$ are \emph{equivalent}, written~$C \sim_I C'$, if there exists a continuous path~$\gamma \colon [0,1] \to \mathcal{S}$ such that
\begin{equation}
\label{eq:contour_equiv}
\begin{aligned}
\gamma(0) &\in C,\\
\gamma(1) &\in C',\\
\gamma(t) &\notin \bigcup_{Q \in \mathcal{C}_{\mathrm{sing}}} Q
\quad \text{for all } t \in (0,1).
\end{aligned}
\end{equation}
\end{definition}
The equivalence classes of~$\sim_I$ are the \emph{interior classes} of~$\mathcal{S}$; their collection forms the \emph{interior decomposition}. Equivalently, let
\begin{equation}
\widehat{\mathcal{S}} \;=\; \mathcal{S} \setminus \bigcup_{Q \in \mathcal{C}_{\mathrm{sing}}} Q
\label{eq:widehat_S}
\end{equation}
denote the resonance surface with all singular contours removed; then $C \sim_I C'$ if and only if~$C$ and~$C'$ lie in the same path component of~$\widehat{\mathcal{S}}$. Geometrically, an interior class is a connected region of~$\mathcal{S}$ on which the regular contours form a single continuous family. The decomposition diagrams in panel~(b) of Figures~\ref{bd:fig1}--\ref{bd:fig4} show these interior classes as shaded regions.

An alternative view of the interior decomposition is given by the boundary diagrams in panel~(c) of Figures~\ref{bd:fig1}--\ref{bd:fig4}, which represent regular contours as horizontal lines and record how each contour connects with the boundary~$B$ of~$\mathcal{S}$. Each boundary diagram thus provides a coarser-grained representation that uses only information on~$B$. The associated resonance tongue arrangement, which we refer to as the boundary decomposition of~$\mathcal{S}$, can be read off directly from the boundary diagram. However, the boundary decomposition captures only those regular contours that meet~$B$; contours that do not meet~$B$ remain `visible' only to the interior decomposition. Note that cases~$\mathrm{I}$--$\mathrm{III}$ in Section~\ref{section:casesIandII} are distinguished by having different boundary decompositions. 

\subsection{Resonance structure and resonance transitions}
\label{subsec:res_struct_res_surf}
The interior decomposition records how the regular contours of~$\mathcal{S}_\eta$ are grouped between singular contours, but it does not record which singular contours bound each interior class. This information may be completed by encoding the interior decomposition as a labeled graph~\cite{edelsbrunner2010computational}, with singular contours as vertices, interior classes as edges, and the incidence relation specifying which singular contours bound each class.
\begin{definition}[Labeled graph of~$\mathcal{S}_\eta$]
\label{def:incidence_graph}
Let~$\mathcal{C}_{\mathrm{sing}}$ denote the set of singular contours of~$\mathcal{S}_\eta$, and let~$\mathcal{I}$ be the collection of interior classes. The \emph{labeled graph} of~$\mathcal{S}_\eta$ is
\begin{equation}
G(\mathcal{S}_\eta) \;=\; (V_\eta, E_\eta, \lambda_\eta),
\label{eq:incidence_graph_new}
\end{equation}
with vertex set~$V_\eta = \mathcal{C}_{\mathrm{sing}}$ and edge set~$E_\eta = \mathcal{I}$. An interior class~$e \in E_\eta$ is incident to a singular contour~$v \in V_\eta$ if and only if~$v \subset \overline{e}$, where the closure $\overline{e}$ of $e$ is taken in~$\mathcal{S}_\eta$. Each vertex~$v \in V_\eta$ carries the label~$\lambda_\eta(v) = \rho_\eta(v)$, which is well-defined since~$\rho_\eta$ is constant on every singular contour.
\end{definition}

Since~$D_\eta$ is compact and~$\rho_\eta$ is continuous, the extreme value theorem implies that every interior class attains its maximum and minimum values of~$\rho_\eta$, each on a singular contour bounding the class. Each edge of~$G(\mathcal{S}_\eta)$ is therefore incident to exactly two vertices, so~$G(\mathcal{S}_\eta)$ is indeed a graph. It provides an alternative representation of the decomposition diagram shown in Figures~\ref{bd:fig1}--\ref{bd:fig4}. Both are combinatorial invariants of~$\mathcal{S}_\eta$ that we will use to compare resonance surfaces in the one-parameter collection~$\{\mathcal{S}_\eta\}$ and to detect the \emph{resonance transitions} at which their qualitative organization changes.

\begin{definition}[Equivalence of resonance surfaces]
\label{def:resonance_equivalence}
Let~$\mathcal{S}_{\eta_1}$ and~$\mathcal{S}_{\eta_2}$ be two resonance surfaces with labeled graphs
\[
G(\mathcal{S}_{\eta_1}) = (V_{\eta_1}, E_{\eta_1}, \lambda_{\eta_1}),
\qquad
G(\mathcal{S}_{\eta_2}) = (V_{\eta_2}, E_{\eta_2}, \lambda_{\eta_2}).
\]
We say that~$\mathcal{S}_{\eta_1}$ and~$\mathcal{S}_{\eta_2}$ are \emph{equivalent}, and write~$\mathcal{S}_{\eta_1} \cong \mathcal{S}_{\eta_2}$, if there exists a graph isomorphism~$\Phi \colon G(\mathcal{S}_{\eta_1}) \to G(\mathcal{S}_{\eta_2})$ such that
\[
\lambda_{\eta_2}(\Phi(v)) < \lambda_{\eta_2}(\Phi(w))
\quad\text{whenever}\quad
\lambda_{\eta_1}(v) < \lambda_{\eta_1}(w),
\]
for all~$v, w \in V_{\eta_1}$.
\end{definition}

This equivalence identifies two resonance surfaces if their singular contours and interior classes can be matched while preserving the rotation-number ordering between them. A \emph{resonance transition} is a value of~$\eta = \eta^*$ at which this equivalence fails between resonance surfaces on either side.
\begin{definition}[Resonance transition]
\label{def:resonance_transition}
A \emph{resonance transition} occurs at a parameter value~$\eta = \eta^*$ if, for all sufficiently small~$r > 0$,
\begin{equation}
\mathcal{S}_{\eta^* - r} \;\not\cong\; \mathcal{S}_{\eta^* + r}.
\label{eq:resonance_transition}
\end{equation}
\end{definition}
The resonance structure at a value~$\eta_0$ is \emph{structurally stable} if it is preserved up to equivalence on an open~$\eta$-interval containing~$\eta_0$; a resonance transition occurs at the boundary of any such interval. This mirrors the standard viewpoint of bifurcation in dynamical systems theory~\cite{wiggins2003introduction} and in one-parameter Morse--Cerf theory~\cite{cerf1970stratification}.

\begin{figure*}[t!]
  \centering
      \includegraphics[width=\linewidth]{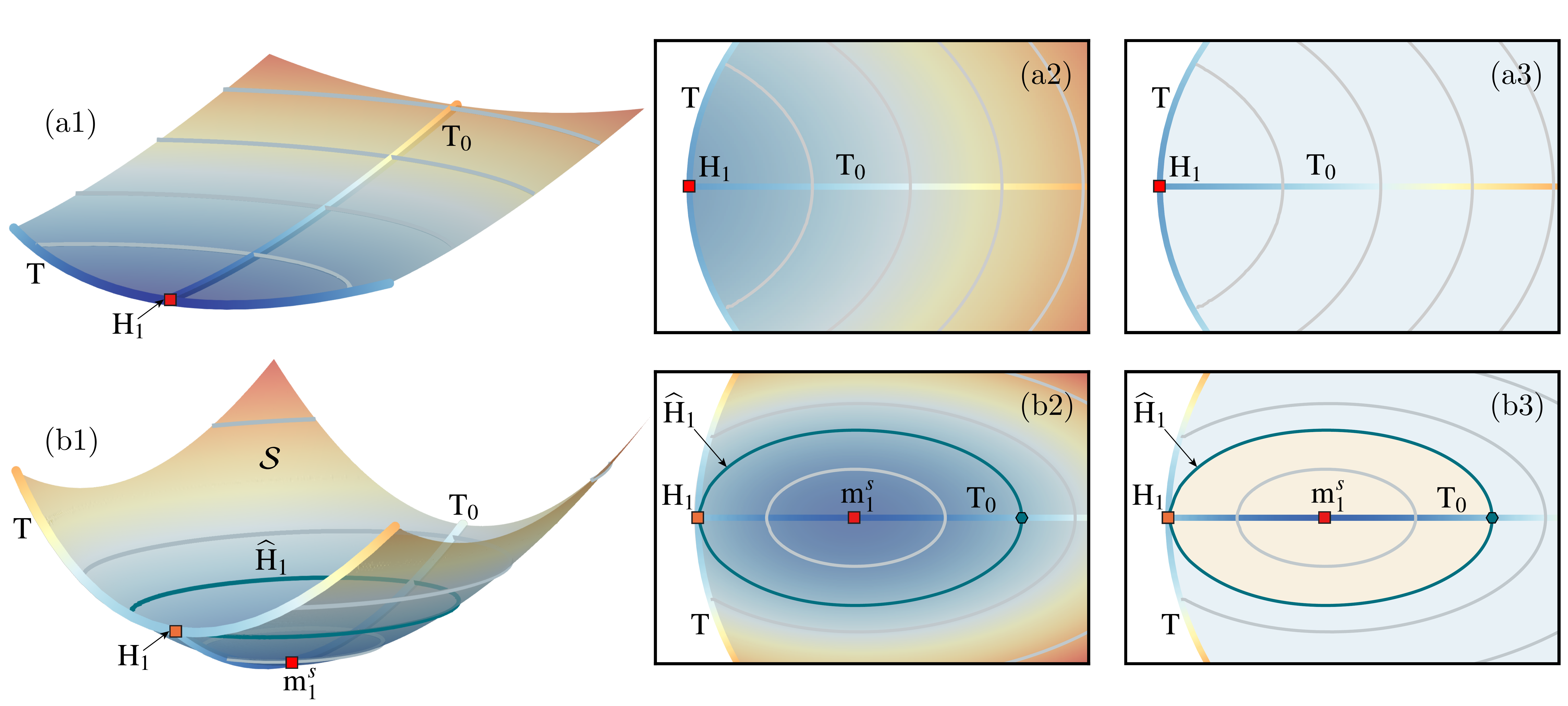}
\caption{\label{fig:H1Trans} The resonance structure for the corner--saddle transition~$\mathrm{HT}_0$ at the corner point~$\mathrm{H}_1$. Panel~(a1), before the transition at $\eta=-0.220$, shows the resonance surface~$\mathcal{S}$, panel~(a2) the resonance diagram, and panel~(a3) the decomposition diagram with a single interior class (blue); $\mathrm{H}_1$ is a corner minimum (red square). Panels~(b1)--(b3), after the transition at $\eta=-0.230$, show~$\mathrm{H}_1$ as a corner saddle (orange square) with separatrix~$\widehat{\mathrm{H}}_1$ (turquoise curve), and the boundary minimum~$\mathrm{m}_1^s$ (red square).}
\end{figure*}

\begin{figure*}[t!]
  \centering
      \includegraphics[width=\linewidth]{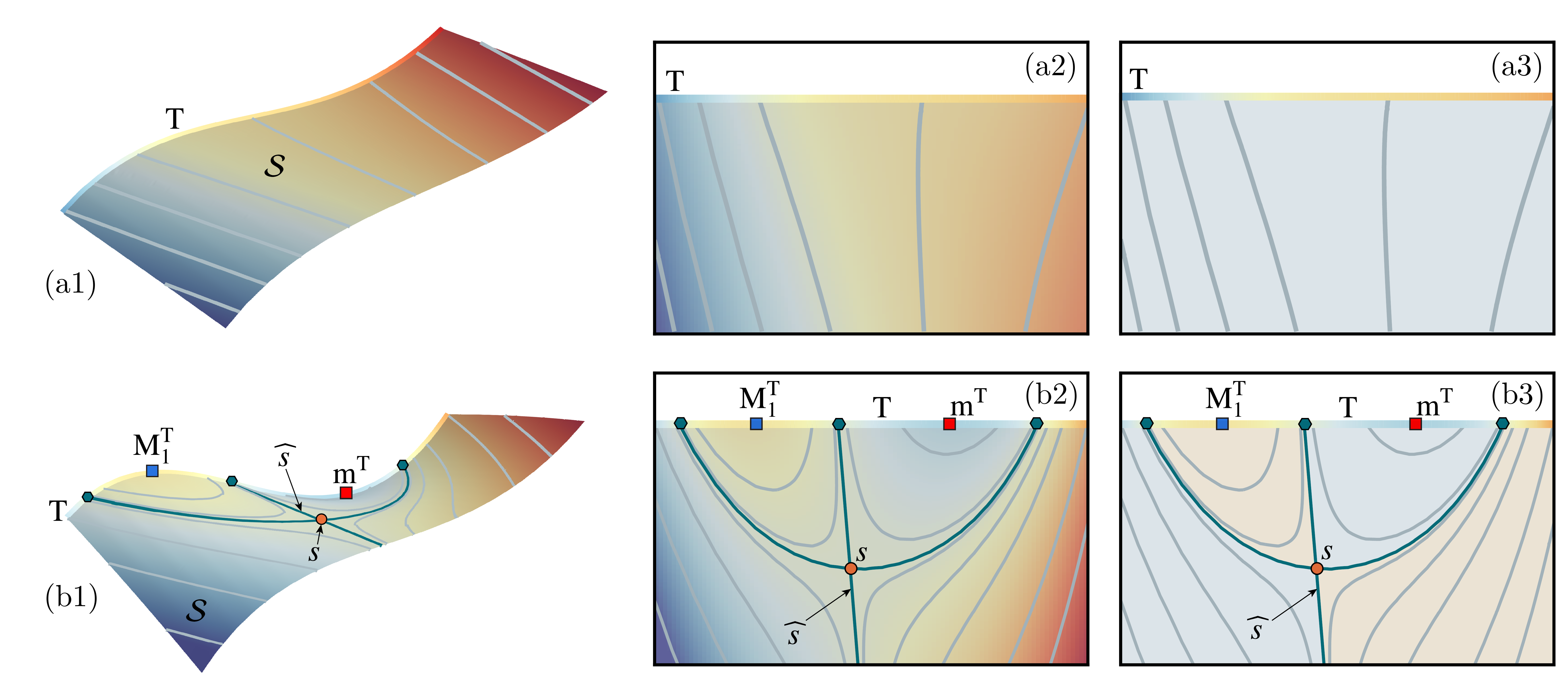}
\caption{\label{fig:CT_trans} The resonance structure for the cusp transition~$\mathrm{CT}$, in the format of Figure~\ref{fig:H1Trans}. Panels~(a1)--(a3), at $\eta=-0.265$, show the resonance structure near~$\mathrm{T}$. Panels~(b1)--(b2), at $\eta=-0.2705$, show the boundary maximum~$\mathrm{M}_1^T$ (blue square) and the boundary minimum~$\mathrm{m}^T$ (red square), together with the interior saddle~$s$ (orange dot) and its separatrix~$\widehat{s}$ (turquoise line); panel~(b3) contains four interior classes shaded in alternating blue and yellow.}
\end{figure*}

\begin{figure*}[t!]
  \centering
      \includegraphics[width=\linewidth]{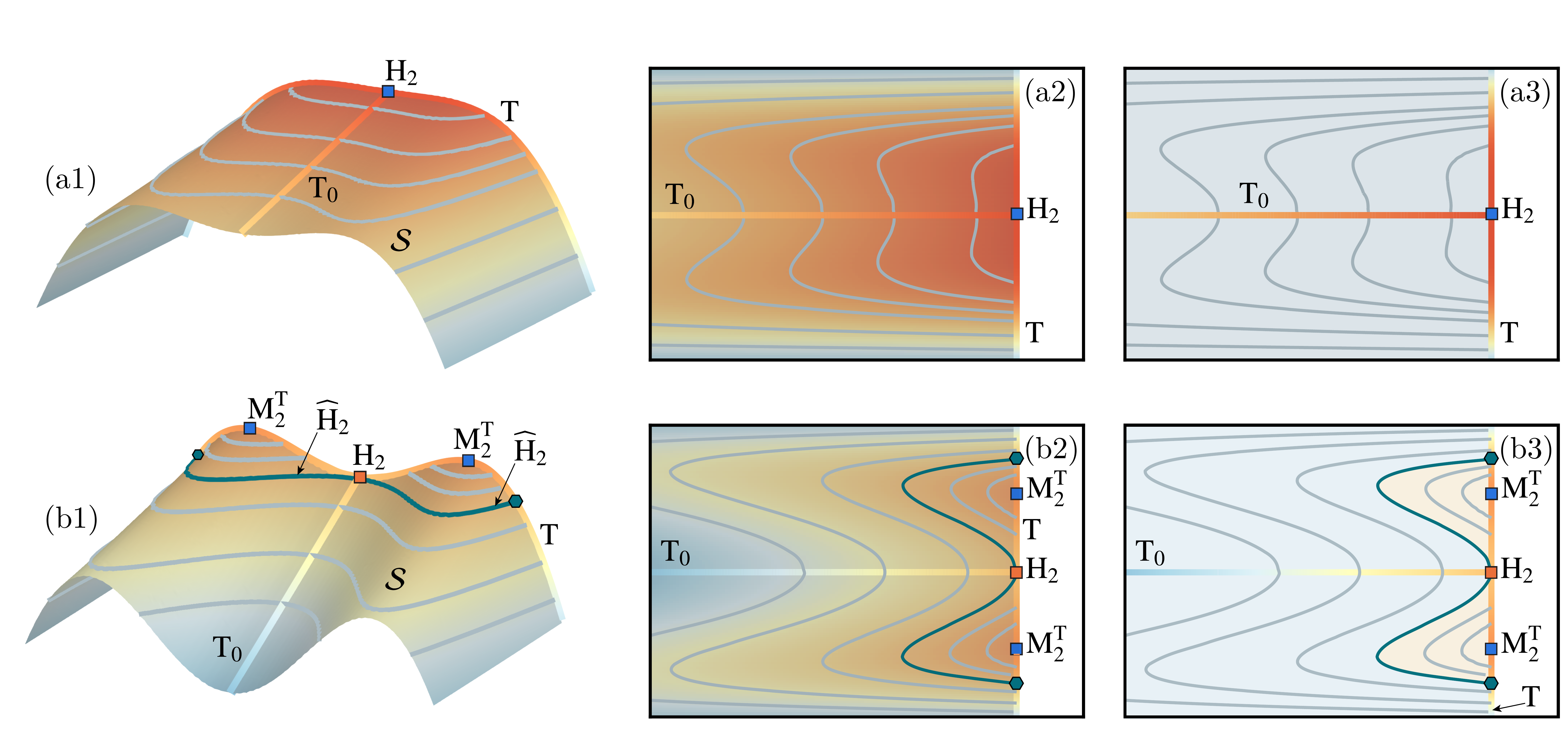}
\caption{\label{fig:H2Trans} The resonance structure for the corner--saddle transition~$\mathrm{HT}$ at the corner point~$\mathrm{H}_2$, in the format of Figure~\ref{fig:CT_trans}. Panels~(a1)--(a3), at $\eta=-0.280$, show the corner maximum~$\mathrm{H}_2$ (blue square) and a single interior class. Panels~(b1)--(b2), at $\eta=-0.2835$, show~$\mathrm{H}_2$ as a corner saddle (orange square) with separatrix~$\widehat{\mathrm{H}}_2$ (turquoise curve), and the boundary maximum~$\mathrm{M}_2^T$ (blue dot) on~$\mathrm{T}$; panel~(b3) contains two interior classes.}
\end{figure*}
\begin{figure*}[t!]
  \centering
      \includegraphics[width=\linewidth]{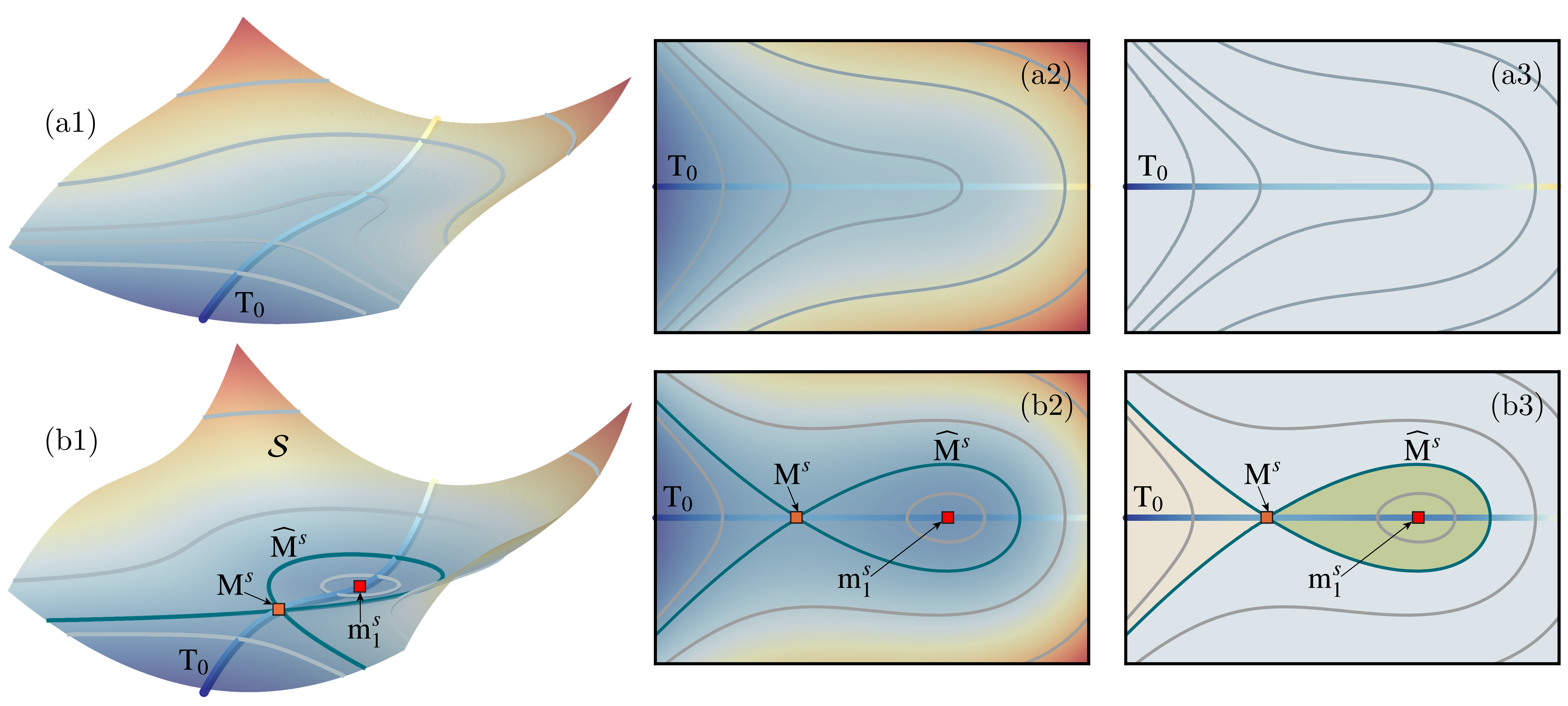}
\caption{\label{fig:transition_saddle_center} The resonance structure for the symmetric cusp transition~$\mathrm{CT}_0$ on the symmetry segment~$\mathrm{T}_0$, in the format of Figure~\ref{fig:CT_trans}. Panels~(a1)--(a3), at $\eta=-0.2794$, show the resonance structure near~$\mathrm{T}_0$.  Panels~(b1)--(b3), at $\eta=-0.2802$, show the boundary saddle~$\mathrm{M}^s$ (orange square) with separatrix~$\widehat{\mathrm{M}}^s$ (turquoise curve), and the boundary minimum~$\mathrm{m}_2^s$ (red square); panel~(b3) has three interior classes.}
\end{figure*}

\section{Resonance transitions in the periodically forced Welander model}
\label{section:res_bd_unfold}
We now consider the resonance transitions of system~\eqref{eq:forced_ODE_second} by comparing the resonance structure at values of~$\eta$ immediately before and after each transition. We use three complementary representations: a local view of the resonance surface~$\mathcal{S}_\eta$ rendered from computed data, the corresponding resonance diagram in the $(\mu,c)$-plane, and the decomposition diagram. Section~\ref{section:res_bd_unfold_1} presents the two transitions between cases~$\mathrm{I}$--$\mathrm{III}$, and Section~\ref{section:res_bd_unfold_2} examines two further transitions, which concern the additional cases~$\mathrm{IV}$--$\mathrm{VI}$ described later in Section~\ref{section:cases_organisation}. Together, these four transitions account for the passage from case~$\mathrm{I}$ to case~$\mathrm{VI}$ at the level of the boundary decomposition.

\subsection{Resonance transitions between cases~$\mathrm{I}$--$\mathrm{III}$}
\label{section:res_bd_unfold_1}
The two transitions between cases~$\mathrm{I}$--$\mathrm{III}$ are a \emph{corner--saddle transition}~$\mathrm{HT}_0$ between cases~$\mathrm{I}$ and~$\mathrm{II}$, which creates a boundary minimum on~$\mathrm{T}_0$; and a \emph{cusp transition}~$\mathrm{CT}$ between cases~$\mathrm{II}$ and~$\mathrm{III}$, which creates a pair of boundary extrema on the torus bifurcation curve~$\mathrm{T}$.

Figure~\ref{fig:H1Trans} shows the corner--saddle transition~$\mathrm{HT}_0$ at the corner point~$\mathrm{H}_1$, illustrated for both signs of the forcing amplitude~$c$. Before the transition, panels~(a1)--(a3) at $\eta = -0.220$ show a single family of resonance tongues near the corner minimum~$\mathrm{H}_1$, each connecting a point on~$\mathrm{T}_0$ to a point on~$\mathrm{T}$; this situation is topologically as in Figure~\ref{fig:topo_sketches_H}(a). After the transition, panels~(b1)--(b3) at $\eta = -0.230$ show that~$\mathrm{H}_1$ has become a corner saddle, and that a boundary minimum~$\mathrm{m}_1^s$ has emerged from~$\mathrm{H}_1$ along~$\mathrm{T}_0$. Panel~(b3) in particular shows that the separatrix~$\widehat{\mathrm{H}}_1$ separates the original family of resonance tongues connecting~$\mathrm{T}_0$ to~$\mathrm{T}$ from a new family that is bounded by~$\widehat{\mathrm{H}}_1$ and organized by~$\mathrm{m}_1^s$. Note that this corner--saddle transition is responsible for the difference in the boundary decompositions of Figures~\ref{bd:fig1} and~\ref{bd:fig2}.

\begin{figure*}[t!]
  \centering
  \includegraphics[width=\linewidth]{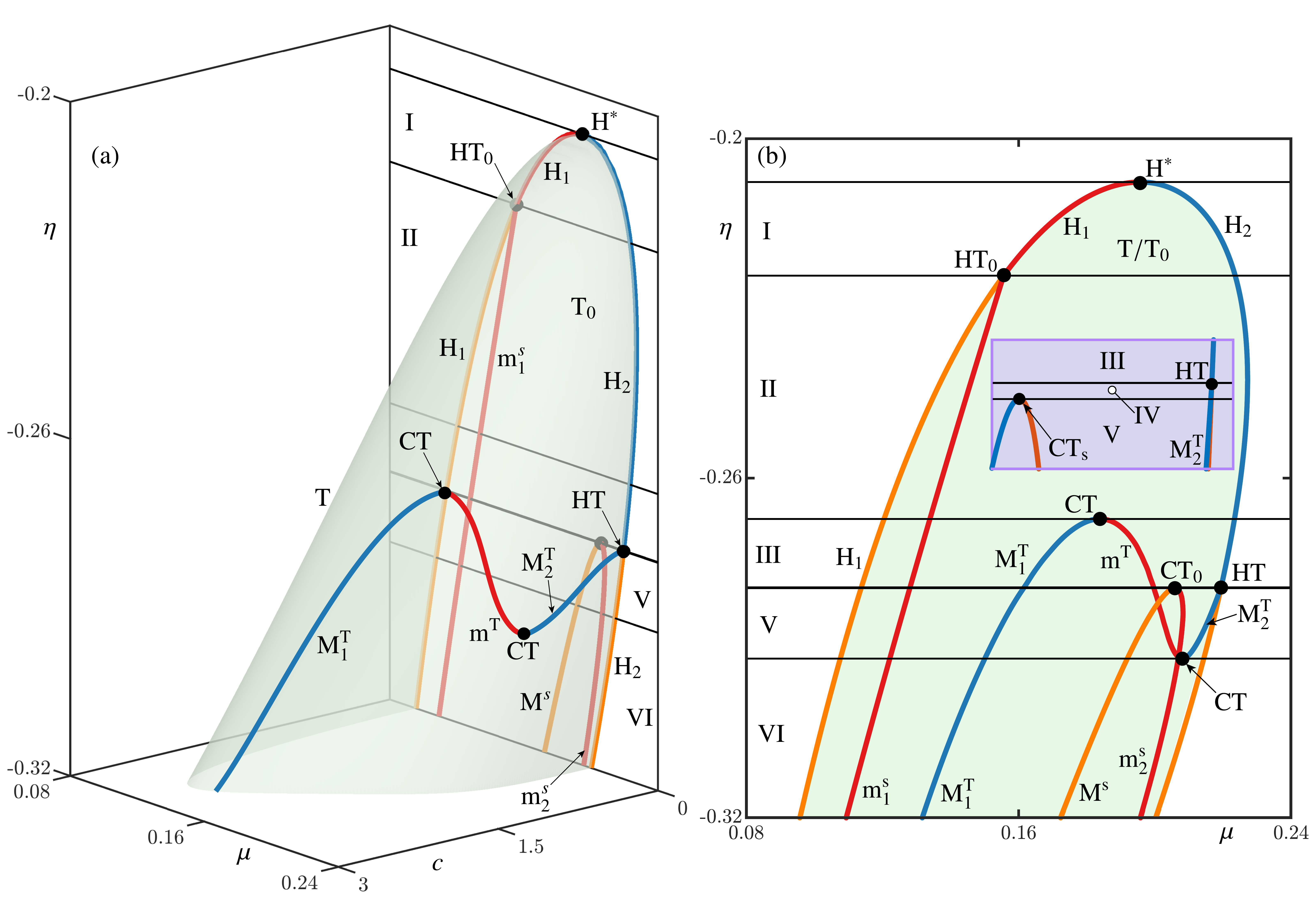}
\caption{\label{fig:refinedBD}
Curves of boundary and corner singularities in $(\mu,\eta,c)$-space in panel~(a), and their projection onto the $(\mu,\eta)$-plane in panel~(b). Shown are the torus bifurcation surface~$\mathrm{T}$ and the symmetry set~$\mathrm{T}_0$ (shaded green). Also shown are the boundary-maximum curves~$\mathrm{M}_1^{T}$ and~$\mathrm{M}_2^{T}$, and the boundary-minimum curve~$\mathrm{m}^{T}$ on~$\mathrm{T}$. The boundary-saddle curve~$\mathrm{M}^{s}$ and the boundary-minimum curves~$\mathrm{m}_1^{s}$ and~$\mathrm{m}_2^{s}$ are shown on~$\mathrm{T}_0$, and the curves~$\mathrm{H}_1$ and~$\mathrm{H}_2$ in the intersection of~$\mathrm{T}$ and~$\mathrm{T}_0$. These curves are colored blue for maxima, red for minima, and orange for saddles, and meet at the points~$\mathrm{H}^{*}$,~$\mathrm{HT}_0$,~$\mathrm{HT}$,~$\mathrm{CT}$, and~$\mathrm{CT}_0$. The corresponding resonance transitions are indicated by black horizontal lines, which partition the $\eta$-axis into cases~$\mathrm{I}$--$\mathrm{VI}$. Panel~(b) includes an inset enlarging the narrow parameter region containing case~$\mathrm{IV}$.}
\end{figure*}

Figure~\ref{fig:CT_trans} illustrates the cusp transition~$\mathrm{CT}$, rendered in curvilinear coordinates along the relevant portion of~$\mathrm{T}$. Before the transition, panels~(a1)--(a3) at $\eta = -0.265$ show that the region of~$\mathcal{S}$ near~$\mathrm{T}$ is free of singularities: contours either end on~$\mathrm{T}$ or have no endpoints there. After the transition, panels~(b1)--(b3) at $\eta = -0.2705$ show that a boundary maximum~$\mathrm{M}_1^T$ and a boundary minimum~$\mathrm{m}^T$ have emerged on~$\mathrm{T}$; this happens at a cusp (cubic) singularity~\cite{arnold1981singularity} of the rotation number~$\rho$ restricted to~$\mathrm{T}$. An important ingredient of the transition is the simultaneous emergence of a $(2,2)$-compound saddle~$s$ in the interior of~$D$. Panel~(b1) shows that~$s$ and its separatrix~$\widehat{s}$ together divide~$\mathcal{S}$ into two components above and two below; panels~(b2)--(b3) show the four corresponding families of resonance tongues, including the two new families organized by~$\mathrm{M}_1^T$ and~$\mathrm{m}^T$.

The resonance transitions~$\mathrm{HT}_0$ and~$\mathrm{CT}$ are distinguished by changes in the maxima and minima on the boundary~$B$; we refer to such transitions as \emph{boundary transitions}.

\subsection{Additional boundary transitions}
\label{section:res_bd_unfold_2}
In direct analogy to the two boundary transitions of Section~\ref{section:res_bd_unfold_1}, we consider two further ones, with the roles of~$\mathrm{T}$ and~$\mathrm{T}_0$ interchanged: a corner--saddle transition~$\mathrm{HT}$, which creates a boundary maximum on~$\mathrm{T}$; and a symmetric cusp transition~$\mathrm{CT}_0$, which creates a boundary saddle and a boundary minimum together on~$\mathrm{T}_0$.

Figure~\ref{fig:H2Trans} illustrates the transition~$\mathrm{HT}$ in curvilinear coordinates along the relevant portion of~$\mathrm{T}$. Before the transition, panels~(a1)--(a3) at $\eta = -0.280$ show a single family of resonance tongues encircling the corner maximum~$\mathrm{H}_2$, analogous to the family around the corner minimum~$\mathrm{H}_1$ in Figure~\ref{fig:H1Trans}(a1)--(a3). After the transition, panels~(b1)--(b3) at $\eta = -0.2835$ show that~$\mathrm{H}_2$ has become a corner saddle and that a boundary maximum~$\mathrm{M}_2^T$ has emerged from~$\mathrm{H}_2$ along~$\mathrm{T}$, with one copy on each side of~$\mathrm{T}_0$. Panel~(b3) shows that the separatrix~$\widehat{\mathrm{H}}_2$ has two branches, each bounding a new family of resonance tongues organized around a copy of~$\mathrm{M}_2^T$.

Figure~\ref{fig:transition_saddle_center} presents the transition~$\mathrm{CT}_0$ near~$\mathrm{T}_0$. Before the transition, panels~(a1)--(a3) at $\eta = -0.2794$ show a single family of resonance tongues, resembling that in Figure~\ref{fig:CT_trans}(a1)--(a3) but extended across~$\mathrm{T}_0$ by reflection. After the transition, panels~(b1)--(b3) at $\eta = -0.2802$ show that a boundary saddle~$\mathrm{M}^s$ and a boundary minimum~$\mathrm{m}_2^s$ have emerged on~$\mathrm{T}_0$ via a cusp singularity of~$\rho$ restricted to~$\mathrm{T}_0$. Unlike the cusp transition on~$\mathrm{T}$ in Figure~\ref{fig:CT_trans},~$\mathrm{CT}_0$ is constrained by reflection symmetry, so the new singularities are confined to~$\mathrm{T}_0$ and no interior singularity is created. Panel~(b3) shows that the separatrix~$\widehat{\mathrm{M}}^s$ separates three families of resonance tongues: two with a single endpoint on~$\mathrm{T}_0$ that lie respectively above and below~$\widehat{\mathrm{M}}^s$ on~$\mathcal{S}$, and a third encircling~$\mathrm{m}_2^s$, with both endpoints on~$\mathrm{T}_0$, one on either side of~$\mathrm{m}_2^s$.

\section{Tracking the boundary transitions in $(\mu,\eta,c)$-space}
\label{section:res_bd}
We now locate the boundary transitions of Section~\ref{section:res_bd_unfold} by tracking the boundary singularities as curves. They are shown in the full three-dimensional~$(\mu,\eta,c)$-space in Figure~\ref{fig:refinedBD}(a) and in projection onto the~$(\mu,\eta)$-plane in Figure~\ref{fig:refinedBD}(b). In~$(\mu,\eta,c)$-space, the torus bifurcation~$\mathrm{T}$ is a surface in $c > 0$ and, similarly,~$\mathrm{T}_0$ is a region of the~$(\mu,\eta)$-plane. The two surfaces~$\mathrm{T}$ and~$\mathrm{T}_0$ meet along the corner point curves~$\mathrm{H}_1$ and~$\mathrm{H}_2$, and the boundary extrema trace out curves on them. We locate the boundary extrema for fixed~$\eta$, which is an easy task: the rotation number~$\rho$ varies smoothly on the surfaces~$\mathrm{T}$ and~$\mathrm{T}_0$; it fails to be smooth only along the corner point curves~$\mathrm{H}_1$ and~$\mathrm{H}_2$; see Figure~\ref{fig:intro_theory_fig}(c). We then follow the different extrema as~$\eta$ varies and locate the boundary transitions where these curves meet at folds in~$\eta$, or where they end on~$\mathrm{H}_1$ or~$\mathrm{H}_2$. Because it relies only on data along~$B$, this procedure is not specific to the periodically forced Welander model but applies to any periodically forced system with a region of normally hyperbolic invariant tori.

The curves~$\mathrm{H}_1$ and~$\mathrm{H}_2$ meet at the point~$\mathrm{H}^*$ at~$c=0$ and jointly bound the region~$\mathrm{T}_0$ in Figure~\ref{fig:refinedBD}(b); note that this is also the projection onto the~$(\mu,\eta)$-plane of the surface~$\mathrm{T}$ shown in panel~(a). Each of the curves~$\mathrm{H}_1$ and~$\mathrm{H}_2$ is divided into a segment of corner extrema and a segment of corner saddles by the corner--saddle transitions~$\mathrm{HT}_0$ on~$\mathrm{H}_1$ and~$\mathrm{HT}$ on~$\mathrm{H}_2$, respectively.

On~$\mathrm{T}_0$, the curves~$\mathrm{m}_2^{s}$ and~$\mathrm{M}^{s}$ join to form a single parabolic curve whose maximum in~$\eta$ is the symmetric cusp transition~$\mathrm{CT}_0$; the curve~$\mathrm{m}_1^{s}$ terminates on~$\mathrm{H}_1$ at the corner--saddle transition~$\mathrm{HT}_0$. The surface~$\mathrm{T}$ emerges from~$\mathrm{H}_1\cup\mathrm{H}_2$, and its boundary extrema form a single curve whose two folds in~$\eta$, which are the cusp transitions~$\mathrm{CT}$, divide this curve of extrema into the segments~$\mathrm{M}_1^{T}$,~$\mathrm{m}^{T}$, and~$\mathrm{M}_2^{T}$; the latter branch terminates on~$\mathrm{H}_2$ at the corner--saddle transition~$\mathrm{HT}$.

Cases~$\mathrm{I}$--$\mathrm{III}$ exist as intervals of~$\eta$, separated by the boundary transitions~$\mathrm{HT}_0$ and~$\mathrm{CT}$ from Section~\ref{section:res_bd_unfold_1}. We find three additional cases~$\mathrm{IV}$--$\mathrm{VI}$, separated by the boundary transitions~$\mathrm{HT}$ and~$\mathrm{CT}_0$ from Section~\ref{section:res_bd_unfold_2}, as well as a second cusp transition, also labeled~$\mathrm{CT}$. Note that case~$\mathrm{IV}$ occupies a very narrow~$\eta$-interval, which is enlarged in the inset of Figure~\ref{fig:refinedBD}(b). All five resonance transitions and their~$\eta$-values are listed in Table~\ref{table:res_trans_table}. Thus, Figure~\ref{fig:refinedBD} provides a coarse-grained picture of the resonance structure at the level of the boundary decomposition; the organization of the interior singularities is not captured here. There are numerous interior transitions, whose discussion is left for future work.

\begin{table}[t!]
\centering
\caption{\label{table:res_trans_table}The point~$\mathrm{H}^{*}$ and the five resonance transitions between cases~$\mathrm{I}$--$\mathrm{VI}$.}
\setlength{\tabcolsep}{0.75em}
\begin{tabular}{llcr}
\hline\hline
label & cases & type & $\eta$-value \\
\hline
$\mathrm{H}^{*}$  & $\emptyset$--$\mathrm{I}$      & --             & $-0.207676$ \\
$\mathrm{HT}_0$ & $\mathrm{I}$--$\mathrm{II}$    & corner--saddle & $-0.224224$ \\
$\mathrm{CT}$     & $\mathrm{II}$--$\mathrm{III}$  & cusp           & $-0.267215$ \\
$\mathrm{HT}$   & $\mathrm{III}$--$\mathrm{IV}$  & corner--saddle & $-0.279324$ \\
$\mathrm{CT}_0$   & $\mathrm{IV}$--$\mathrm{V}$    & symmetric cusp & $-0.279458$ \\
$\mathrm{CT}$     & $\mathrm{V}$--$\mathrm{VI}$    & cusp           & $-0.291904$ \\
\hline\hline
\end{tabular}
\end{table}

\section{Cases~$\mathrm{III}$--$\mathrm{VI}$ of resonance structure}
\label{section:cases_organisation}
We now describe the resonance structure of cases~$\mathrm{IV}$--$\mathrm{VI}$ identified in Figure~\ref{fig:refinedBD}. To this end, we begin with case~$\mathrm{III}$ at a value of~$\eta$ near~$\mathrm{HT}$, different from that in Section~\ref{section:casesIandII}. Each case is presented again by the resonance diagram in the~$(\mu,c)$-plane, the decomposition diagram, and the boundary diagram.

\begin{figure*}[ht!]
  \centering
  \includegraphics{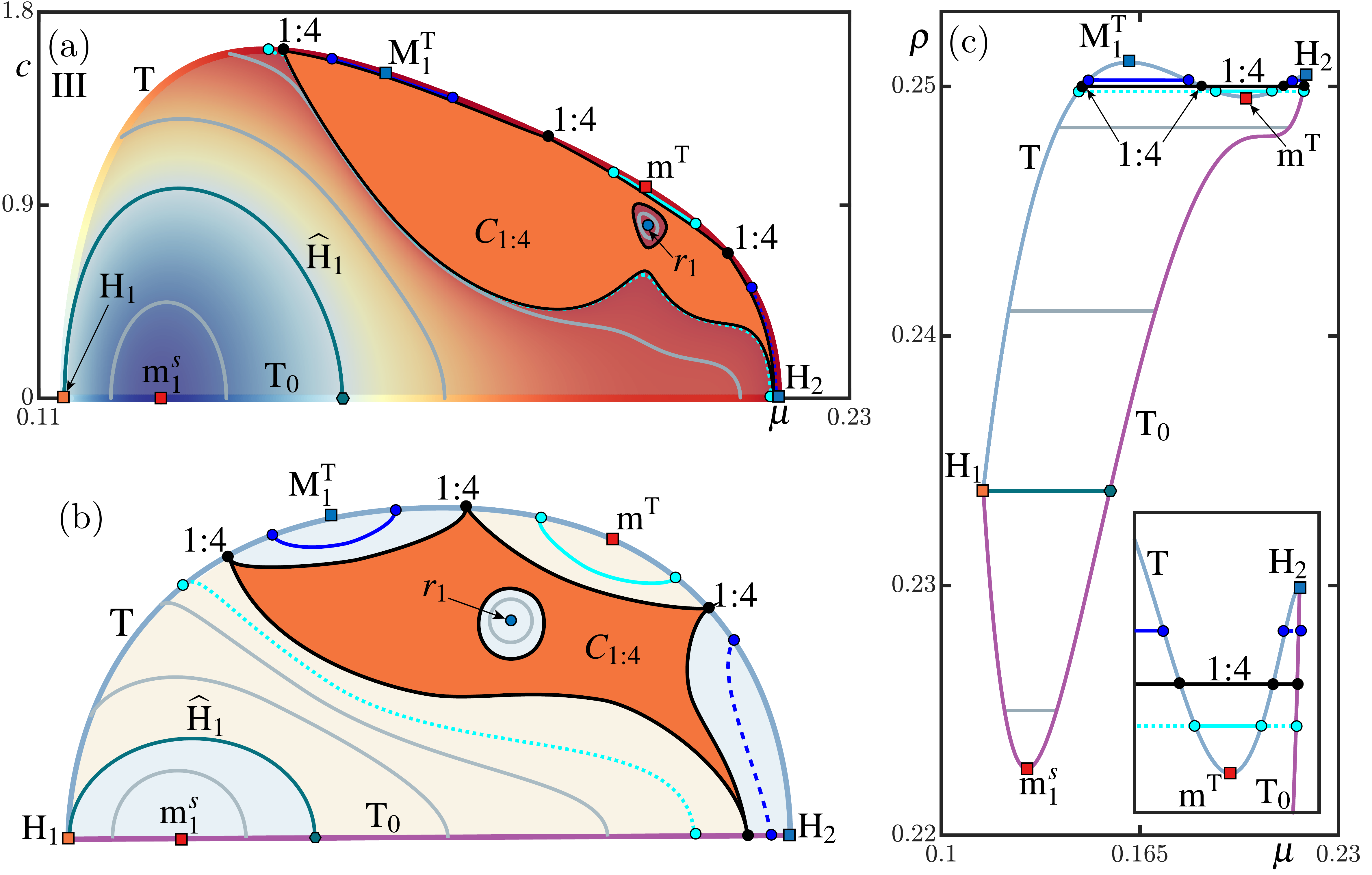}
\caption{The resonance structure for case~$\mathrm{III}$ at $\eta=-0.2793$, in the format of Figures~\ref{bd:fig1}--\ref{bd:fig4}, featuring the prominent~$1{:}4$ resonance tongue~$C_{1{:}4}$. Selected pairs of resonance tongues are shown in cyan and blue, each with a solid and a dotted branch. Panel~(b) shows five families of resonance tongues (shaded in blue and yellow). Panel~(c) contains an inset that enlarges a region near~$\mathrm{m}^{T}$; the value $\rho = 1/4$ is indicated by the black line.}
\label{bd:fig3}
\end{figure*}

\subsection{Case~$\mathrm{III}$ revisited}
\label{section:res_case_III}
Figure~\ref{bd:fig3} shows case~$\mathrm{III}$ at the lower value of $\eta = -0.2793$, close to the boundary transition~$\mathrm{HT}$. The boundary singularities are as in Figure~\ref{bd:fig4}; they are the boundary maximum~$\mathrm{M}_1^T$; the boundary minima~$\mathrm{m}^T$ and~$\mathrm{m}_1^s$; the corner maximum~$\mathrm{H}_2$; and the corner saddle~$\mathrm{H}_1$. The arrangements of resonance tongues near each of these points are also unchanged; compare Figures~\ref{bd:fig4}(b) and~\ref{bd:fig3}(b).

What is different is that the saddle-node curves of the~$1{:}4$ resonance tongue~$C_{1{:}4}$ now meet~$\mathrm{T}$ at three points and~$\mathrm{T}_0$ at one point, delimiting the regions around~$\mathrm{M}_1^T$,~$\mathrm{m}^T$, and~$\mathrm{H}_2$. In doing so,~$C_{1{:}4}$ plays the role of the separatrix~$\widehat{s}_1$ in Figure~\ref{bd:fig4}, separating the same five families of resonance tongues around the boundary singularities on~$B$. An additional interior maximum~$r_1$ has appeared, which lies in a circular region of larger~$\rho$ `inside'~$C_{1{:}4}$; note, however, that the existence of this `local' region does not impede the separating role of~$C_{1{:}4}$. On the associated resonance surface (not shown), three of these families lie above~$C_{1{:}4}$ (including the family around~$r_1$) and two below. The contour~$C_{1{:}4}$ is therefore a compound saddle of order~$(2,3)$ according to Definition~\ref{def:local_singularity_type}(iv).

The resonance structure is otherwise the same, as shown by the boundary diagram in Figure~\ref{bd:fig3}(c). Below the line at~$\rho=1/4$, the overlapping cyan pair corresponds to the family around~$\mathrm{m}^T$, with both endpoints on~$\mathrm{T}$, and a family connecting~$\mathrm{T}_0$ to~$\mathrm{T}$. Above this line, the disjoint blue pair corresponds to the family around~$\mathrm{M}_1^T$, also with both endpoints on~$\mathrm{T}$, and the family around~$\mathrm{H}_2$, connecting~$\mathrm{T}_0$ to~$\mathrm{T}$.

\begin{figure*}[t!]
  \centering
  \includegraphics{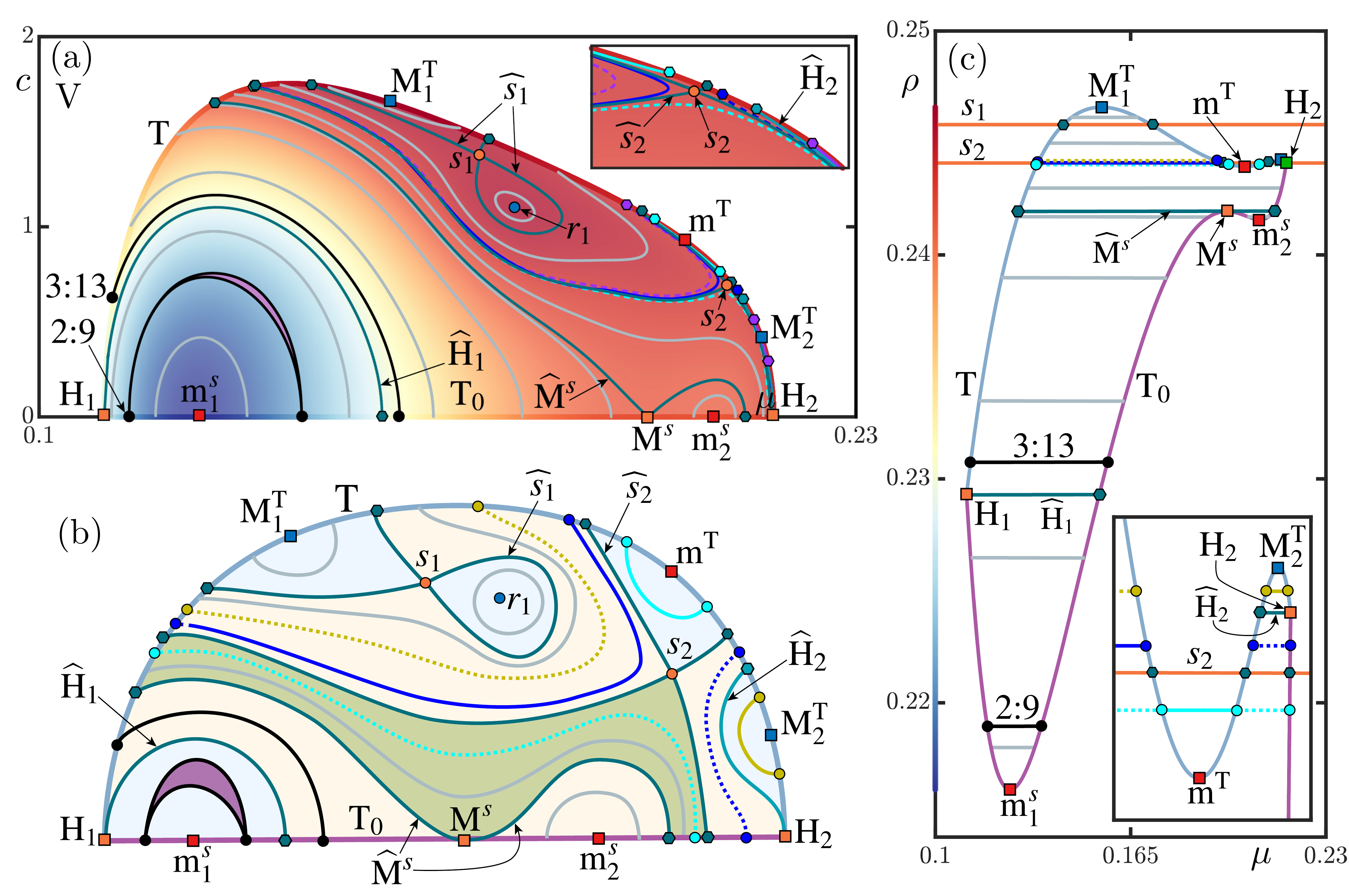}
\caption{\label{bd:caseIVV} The resonance structure for cases~$\mathrm{IV}$ and~$\mathrm{V}$ at $\eta=-0.286$, in the format of Figure~\ref{bd:fig3}. Additionally shown are the maximum~$\mathrm{M}_2^{T}$ (blue square), the boundary saddle~$\mathrm{M}^{s}$ (orange square); the boundary minimum~$\mathrm{m}_2^{s}$ (red square); the interior saddles~$s_1$ and~$s_2$ (orange circles); and the corner saddle~$\mathrm{H}_2$ (orange dot). Shown also are the separatrices~$\widehat{\mathrm{M}}^{s}$,~$\widehat{s}_1$ and~$\widehat{s}_2$, and~$\widehat{\mathrm{H}}_2$ (turquoise). Three selected pairs of resonance tongues at distinct values of~$\rho$ are shown in blue, cyan, and gold, each with a solid and a dotted branch; also shown are the $2{:}9$ and $3{:}13$ resonance tongues. Panel~(b) now shows ten families of resonance tongues (shaded in yellow, blue, and green). Panel~(c) contains an inset that enlarges a region near~$\mathrm{m}^{T}$ and~$\mathrm{M}_2^{T}$; the $\rho$-values of~$s_1$ and~$s_2$ are indicated by the orange lines.}
\end{figure*}

\begin{figure*}[t!]
  \centering
  \includegraphics{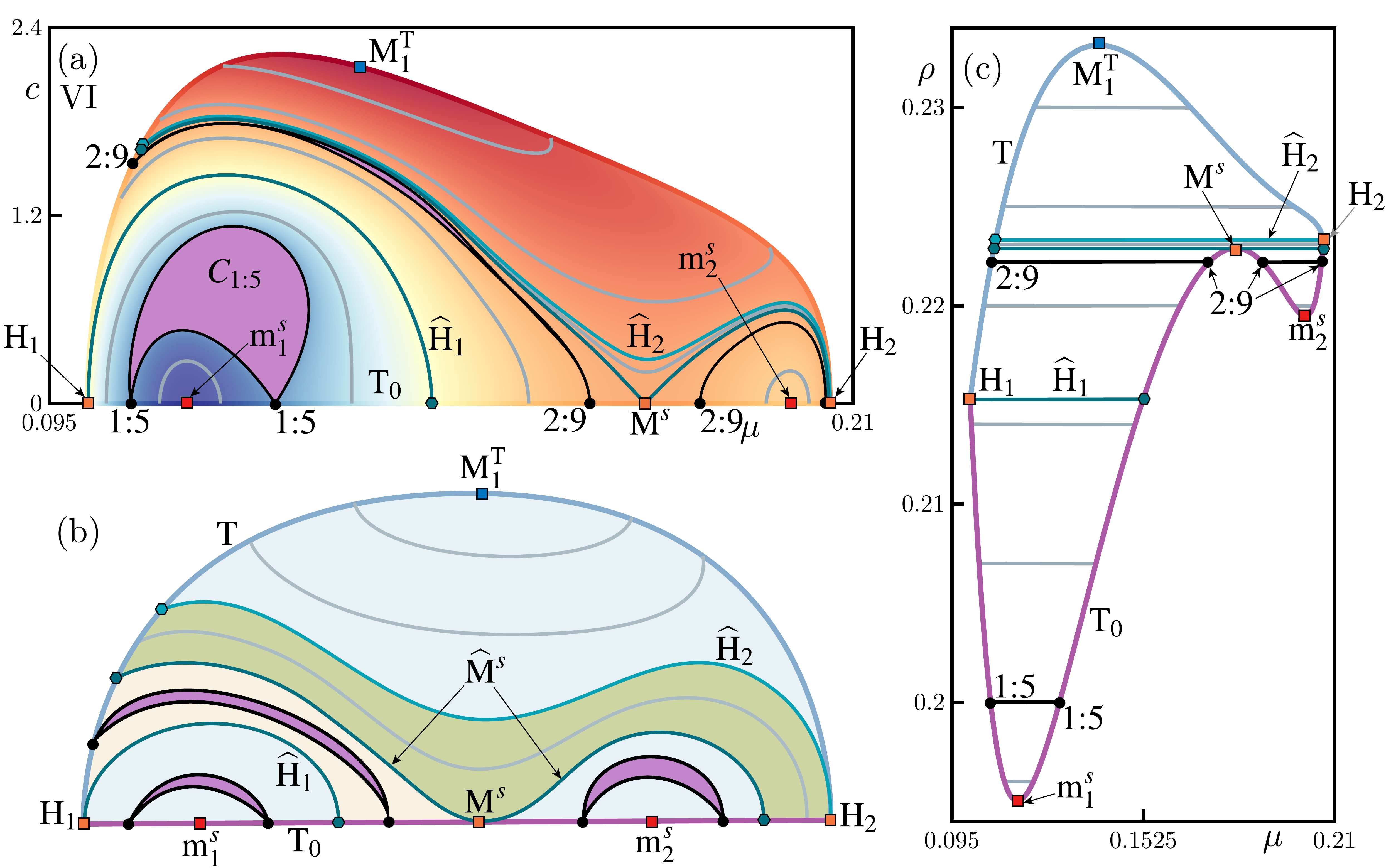}
\caption{\label{bd:caseIVVb} The resonance structure for case~$\mathrm{VI}$ at $\eta=-0.308$, in the format of Figure~\ref{bd:caseIVV}. Panels~(a)--(c) feature the $1{:}5$ and $2{:}9$ resonance tongues; panel~(b) shows five interior classes.}
\end{figure*}

\subsection{Cases~$\mathrm{IV}$ and~$\mathrm{V}$}
Since case~$\mathrm{IV}$ occupies only a narrow $\eta$-interval and is therefore difficult to resolve, we illustrate the resulting resonance structure of cases~$\mathrm{IV}$ and~$\mathrm{V}$ together in Figure~\ref{bd:caseIVV} at $\eta = -0.286$.

What is new along the boundary~$B$ in case~$\mathrm{IV}$ is the appearance of the boundary maximum~$\mathrm{M}_2^T$ on~$\mathrm{T}$ and the change of~$\mathrm{H}_2$ to a corner saddle, following the corner--saddle transition~$\mathrm{HT}$ in Figure~\ref{fig:H2Trans}. As a result, a new family of resonance tongues appears in Figure~\ref{bd:caseIVV}(b), delimited by the separatrix~$\widehat{\mathrm{H}}_2$ and represented by the solid-gold line in Figure~\ref{bd:caseIVV}(c). Case~$\mathrm{V}$ retains these features and additionally developed a boundary saddle~$\mathrm{M}^s$ and a boundary minimum~$\mathrm{m}_2^s$ on~$\mathrm{T}_0$, following the symmetric cusp transition~$\mathrm{CT}_0$ in Figure~\ref{fig:transition_saddle_center}. The separatrix~$\widehat{\mathrm{M}}^s$ delimits the family of resonance tongues forming around~$\mathrm{m}_2^s$ and separates two further families that connect~$\mathrm{T}_0$ to~$\mathrm{T}$.

The interior of~$D$ is now divided by the separatrices~$\widehat{s}_1$ and~$\widehat{s}_2$ of a monotone saddle~$s_1$ and a~$(2,2)$-compound saddle~$s_2$. These two saddles replace the~$(2,3)$-compound saddle of case~$\mathrm{III}$ and appear as horizontal lines at their respective~$\rho$-values in Figure~\ref{bd:caseIVV}(c). The separatrix~$\widehat{s}_1$ encloses the families around~$\mathrm{M}_1^T$ and~$r_1$, whereas~$\widehat{s}_2$ encloses the family around~$\mathrm{m}^T$. The two separatrices also bound new regions:~$\widehat{s}_2$ encloses a narrow region containing~$\mathrm{H}_2$ and~$\mathrm{M}_2^T$ in the lower right of~$D$, while~$\widehat{s}_1$ and~$\widehat{s}_2$ together bound a larger region between them, in which resonance tongues connect~$\mathrm{T}_0$ to~$\mathrm{T}$ directly. Together, the separatrices partition the half-disk~$D$ into the ten families of resonance tongues shown in Figure~\ref{bd:caseIVV}(b).

\subsection{Case~$\mathrm{VI}$}
In case~$\mathrm{VI}$, shown in Figure~\ref{fig:refinedBD} for~$\eta = -0.308$, the two boundary extrema~$\mathrm{M}_2^T$ and~$\mathrm{m}^T$ have vanished together with the interior saddle~$s_2$ at the cusp transition~$\mathrm{CT}$. In fact, there are no interior singularities at this~$\eta$-value: the interior saddle~$s_1$ and the interior maximum~$r_1$ disappeared independently, in an interior resonance transition close to~$\mathrm{CT}$. Additionally, just before~$\mathrm{CT}$, an interior heteroclinic connection between the separatrices~$\widehat{\mathrm{H}}_2$ and~$\widehat{s}_2$ in Figure~\ref{bd:caseIVV} causes them to switch places. The further discussion of these interior transitions lies beyond the scope of this paper.

The resonance structure of case~$\mathrm{VI}$ is determined entirely by the remaining boundary singularities, shown in Figure~\ref{bd:caseIVVb}(a)--(b). There are five resonance tongue families, distinguished clearly by how they attach to~$\mathrm{T}_0$ and~$\mathrm{T}$. Two have both endpoints on~$\mathrm{T}_0$ and form around the boundary minima~$\mathrm{m}_1^s$ and~$\mathrm{m}_2^s$; these contain, respectively, the prominent~$1{:}5$ tongue and one of the two thin~$2{:}9$ tongues. Two further families connect~$\mathrm{T}_0$ to~$\mathrm{T}$: one fills the region between the separatrices~$\widehat{\mathrm{H}}_1$ and~$\widehat{\mathrm{M}}^s$ and contains the second~$2{:}9$ tongue; the other fills the narrow region between~$\widehat{\mathrm{M}}^s$ and~$\widehat{\mathrm{H}}_2$ and extends across the half-disk~$D$. The fifth family forms around the boundary maximum~$\mathrm{M}_1^T$, bounded below by the separatrix~$\widehat{\mathrm{M}}^s$, with both endpoints on~$\mathrm{T}$.

\section{Conclusion}
\label{section:conclusion}
We considered the global arrangement of resonance tongues as a Morse-theoretic problem by examining the level-set topology of the resonance surface~$\mathcal{S}_\eta$. For fixed~$\eta$, the resonance tongues are the terraces of~$\mathcal{S}_\eta$ at rational values of~$\rho$, and their arrangement is determined by singular contours. The framework developed here provides a classification of the qualitatively distinct, structurally stable resonance configurations, whereas previously only individual low-order resonance tongues were computed and considered.

For the periodically forced Welander model~\eqref{eq:forced_ODE_second}, we examined the rearrangement of its resonance tongues in the~$(\mu,c)$-plane under variation of~$\eta$ across five boundary transitions. Locating these transitions in $(\mu,\eta,c)$-space classifies the resonance structure into six structurally stable cases~$\mathrm{I}$--$\mathrm{VI}$, distinguished by the singularities of~$\rho$ on the boundary~$B$ of the half-disk~$D$. Resolving the resonance surface accurately enough to detect these singularities required a rotation-number algorithm that bypasses continuation of locked orbits and computes rational~$\rho = p/q$ directly from the stroboscopic map, reaching denominators as high as~$q = 10^4$.

The focus has been on the half-disk~$D$, its boundary~$B$, and the boundary transitions~$\mathrm{HT}_0$, $\mathrm{CT}$, $\mathrm{HT}$, and~$\mathrm{CT}_0$. The classification we presented constitutes a first, quite natural step in understanding the resonance structure. However, it is complete only at the level of the (coarser) boundary decomposition. Already during the passage through cases~$\mathrm{III}$--$\mathrm{V}$, the resonance surface~$\mathcal{S}_\eta$ develops interior singularities. Interior transitions that create and destroy them remain to be classified. Our investigation already hints at a complicated sequence of transitions, including evidence of a heteroclinic-type transition. For values of~$\eta$ below those in Figure~\ref{fig:refinedBD}, we find, in addition, multi-stage sequences of structurally stable configurations separated by pitchfork-like transitions, at which a boundary extremum changes type and resonance tongues emerge or disconnect. A systematic account of these finer transitions of the resonance structure of the Welander model is the subject of ongoing work.

The Welander model served here as a controlled setting in which to demonstrate the organizing role of the singularities of~$\rho$. The framework itself and our approach, however, are not tied to any one model. The rotation-number algorithm has already been implemented for a delay differential equation model of the ENSO~\cite{keane2018chenciner,keane2017climate}, whose solutions evolve in an infinite-dimensional function space, requiring careful construction of the Poincar\'e map~\cite{bolduc2026resonance,bolduc2026enso}. It could likewise be adapted to other non-planar examples, in particular, to coupled-oscillator systems, which exhibit locked dynamics on a torus in much the same way. This will open the possibility of studying their resonance structure in the same spirit.

\section*{Acknowledgments}
This research was supported by Royal Society Te Ap{\=a}rangi Marden Fund grant \#19-UOA-223.

\bibliographystyle{elsarticle-num}
\bibliography{cas-refs}

@article{bailie2025detailed,
  author = {Bailie, John and Dijkstra, Henk A. and Krauskopf, Bernd},
  title = {A detailed analysis of deep-decoupling/deep-coupling oscillations in the {W}elander model},
  journal = {Chaos: An Interdisciplinary Journal of Nonlinear Science},
  volume = {35},
  number = {7},
  pages = {073126},
  year = {2025}
}

@article{detroux2018experimental,
  author = {Detroux, Thibaut and No{\"e}l, Jean-Philippe and Virgin, Lawrence N. and Kerschen, Ga{\"e}tan},
  title = {Experimental study of isolas in nonlinear systems featuring modal interactions},
  journal = {{PLOS ONE}},
  volume = {13},
  number = {3},
  pages = {e0194452},
  year = {2018}
}

@article{kuether2015nonlinear,
  author = {Kuether, Robert J and Renson, Ludovic and Detroux, Thibaut and Grappasonni, Chiara and Kerschen, Ga{\"e}tan and Allen, Mathew S},
  title = {Nonlinear normal modes, modal interactions and isolated resonance curves},
  journal = {Journal of Sound and Vibration},
  volume = {351},
  pages = {299--310},
  year = {2015}
}

@article{marchionne2018synchronisation,
  author = {Marchionne, Arianna and Ditlevsen, Peter and Wieczorek, Sebastian},
  title = {Synchronisation vs. resonance: Isolated resonances in damped nonlinear oscillators},
  journal = {Physica D: Nonlinear Phenomena},
  volume = {380},
  pages = {8--16},
  year = {2018}
}

@article{habib2018isolated,
  author = {Habib, Giuseppe and Cirillo, Giuseppe I. and Kerschen, Ga{\"e}tan},
  title = {Isolated resonances and nonlinear damping},
  journal = {Nonlinear Dynamics},
  volume = {93},
  number = {3},
  pages = {979--994},
  year = {2018}
}

@article{broer2003geometry,
  author = {Broer, Henk W and Golubitsky, Martin and Vegter, Gert},
  title = {The geometry of resonance tongues: a singularity theory approach},
  journal = {Nonlinearity},
  volume = {16},
  number = {4},
  pages = {1511},
  year = {2003}
}

@article{schilder2007computing,
  author = {Schilder, Frank and Peckham, Bruce B},
  title = {Computing {Arnol'd} tongue scenarios},
  journal = {Journal of Computational Physics},
  volume = {220},
  number = {2},
  pages = {932--951},
  year = {2007}
}

@article{tziperman1995irregularity,
  author = {Tziperman, Eli and Cane, Mark A. and Zebiak, Stephen E.},
  title = {Irregularity and locking to the seasonal cycle in an {ENSO} prediction model as explained by the quasi-periodicity route to chaos},
  journal = {Journal of Atmospheric Sciences},
  volume = {52},
  number = {3},
  pages = {293--306},
  year = {1995}
}

@article{marts2007period,
  author = {Marts, Bradley and Simpson, David J. W. and Hagberg, Aric and Lin, Anna L},
  title = {Period doubling in a periodically forced {Belousov--Zhabotinsky} reaction},
  journal = {Physical Review E},
  volume = {76},
  number = {2},
  pages = {026213},
  year = {2007}
}

@article{lin2004resonance,
  author = {Lin, Anna L and Hagberg, Aric and Meron, Ehud and Swinney, Harry L},
  title = {Resonance tongues and patterns in periodically forced reaction-diffusion systems},
  journal = {Physical Review E},
  volume = {69},
  number = {6},
  pages = {066217},
  year = {2004}
}

@article{krauskopf2014bifurcation,
  author = {Krauskopf, Bernd and Sieber, Jan},
  title = {Bifurcation analysis of delay-induced resonances of the {El Ni{\~n}o} Southern Oscillation},
  journal = {Proceedings of the Royal Society A: Mathematical, Physical and Engineering Sciences},
  volume = {470},
  number = {2169},
  pages = {20140348},
  year = {2014}
}

@article{wang2014influence,
  author = {Wang, Hengtong and Sun, Yongjuan and Li, Yichen and Chen, Yong},
  title = {Influence of autapse on mode-locking structure of a {Hodgkin--Huxley} neuron under sinusoidal stimulus},
  journal = {Journal of Theoretical Biology},
  volume = {358},
  pages = {25--30},
  year = {2014}
}

@article{lee2006bifurcation,
  author = {Lee, Sang-Gui and Kim, Seunghwan},
  title = {Bifurcation analysis of mode-locking structure in a {Hodgkin--Huxley} neuron under sinusoidal current},
  journal = {Physical Review E},
  volume = {73},
  number = {4},
  pages = {041924},
  year = {2006}
}

@article{arnold1965small,
  author = {Arnold, Vladimir I},
  title = {Small denominators. I. Mapping of the circumference onto itself},
  journal = {American Mathematical Society Translations, Series 2},
  volume = {46},
  pages = {213--284},
  year = {1965}
}

@article{arnold1977loss,
  author = {Arnold, Vladimir I},
  title = {Loss of stability of self-oscillation close to resonance and versal deformations of equivariant vector fields},
  journal = {Functional Analysis and its Applications},
  volume = {11},
  number = {2},
  pages = {85--92},
  year = {1977}
}

@article{cessi1994simple,
  author = {Cessi, Paola},
  title = {A simple box model of stochastically forced thermohaline flow},
  journal = {Journal of Physical Oceanography},
  volume = {24},
  number = {9},
  pages = {1911--1920},
  year = {1994}
}

@article{bailie2024bifurcation,
  author = {Bailie, John and Krauskopf, Bernd},
  title = {Bifurcation analysis of a conceptual model for vertical mixing in the {North Atlantic}},
  journal = {Physica D: Nonlinear Phenomena},
  volume = {460},
  pages = {134077},
  year = {2024}
}

@article{simonet1994locking,
  author = {Simonet, J. and Warden, M. and Brun, E.},
  title = {Locking and Arnol'd tongues in an infinite-dimensional system: The nuclear magnetic resonance laser with delayed feedback},
  journal = {Physical Review E},
  volume = {50},
  number = {5},
  pages = {3383},
  year = {1994}
}

@article{farokhniaee2017mode,
  author = {Farokhniaee, AmirAli and Large, Edward W.},
  title = {Mode-locking behavior of {Izhikevich} neurons under periodic external forcing},
  journal = {Physical Review E},
  volume = {95},
  number = {6},
  pages = {062414},
  year = {2017}
}

@misc{bolduc2025seasonal,
  author = {Bolduc-St-Aubin, Samuel and Humphries, Antony R},
  title = {Seasonal-forcing-dominated dynamics of a piecewise-smooth {Ghil--Zaliapin--Thompson} {ENSO} model},
  year = {2025},
  note = {arXiv:2510.27084}
}

@incollection{welander1986thermohaline,
  author = {Welander, Pierre},
  title = {Thermohaline effects in the ocean circulation and related simple models},
  booktitle = {Large-scale transport processes in oceans and atmosphere},
  pages = {163--200},
  publisher = {Springer},
  year = {1986}
}

@article{cessi1996convective,
  author = {Cessi, Paola},
  title = {Convective adjustment and thermohaline excitability},
  journal = {Journal of Physical Oceanography},
  volume = {26},
  number = {4},
  pages = {481--491},
  year = {1996}
}

@article{holliday2015multidecadal,
  author = {Holliday, N. Penny and Cunningham, Stuart A. and Johnson, Clare and Gary, S. F. and Griffiths, Colin and Read, J. F. and Sherwin, Toby},
  title = {Multidecadal variability of potential temperature, salinity, and transport in the eastern subpolar {North Atlantic}},
  journal = {Journal of Geophysical Research: Oceans},
  volume = {120},
  number = {9},
  pages = {5945--5967},
  year = {2015}
}

@article{yashayaev2016recurrent,
  author = {Yashayaev, Igor and Loder, John W},
  title = {Recurrent replenishment of {Labrador Sea Water} and associated decadal-scale variability},
  journal = {Journal of Geophysical Research: Oceans},
  volume = {121},
  number = {11},
  pages = {8095--8114},
  year = {2016}
}

@article{das2017quantitative,
  author = {Das, Suddhasattwa and Saiki, Yoshitaka and Sander, Evelyn and Yorke, James A},
  title = {Quantitative quasiperiodicity},
  journal = {Nonlinearity},
  volume = {30},
  number = {11},
  pages = {4111--4140},
  year = {2017}
}

@article{broer2000resonance,
  author = {Broer, Henk and Sim{\'o}, Carles},
  title = {Resonance tongues in {H}ill's equations: a geometric approach},
  journal = {Journal of Differential Equations},
  volume = {166},
  number = {2},
  pages = {290--327},
  year = {2000}
}

@incollection{takens2001forced,
  author = {Takens, Floris},
  title = {Forced oscillations and bifurcations},
  booktitle = {Global analysis of dynamical systems},
  pages = {1--61},
  publisher = {Institute of Physics Publishing},
  year = {2001}
}

@article{cerf1970stratification,
  author = {Cerf, Jean},
  title = {La stratification naturelle des espaces de fonctions diff{\'e}rentiables r{\'e}elles et le th{\'e}oreme de la pseudo-isotopie},
  journal = {Publications Math{\'e}matiques de l'IH{\'E}S},
  volume = {39},
  pages = {5--173},
  year = {1970}
}

@article{carr2003computing,
  author = {Carr, Hamish and Snoeyink, Jack and Axen, Ulrike},
  title = {Computing contour trees in all dimensions},
  journal = {Computational Geometry},
  volume = {24},
  number = {2},
  pages = {75--94},
  year = {2003}
}

@article{keane2018chenciner,
  author = {Keane, Andrew and Krauskopf, Bernd},
  title = {{Chenciner} bubbles and torus break-up in a periodically forced delay differential equation},
  journal = {Nonlinearity},
  volume = {31},
  number = {6},
  pages = {R165--R187},
  year = {2018}
}

@book{arnold1981singularity,
  author = {Arnold, Vladimir Igorevich},
  title = {Singularity theory},
  volume = {53},
  publisher = {Cambridge University Press},
  year = {1981}
}

@book{edelsbrunner2010computational,
  author = {Edelsbrunner, Herbert and Harer, John},
  title = {Computational Topology: an introduction},
  publisher = {American Mathematical Soc.},
  year = {2010}
}

@book{lee2000introduction,
  author = {Lee, John M},
  title = {Introduction to Topological Manifolds},
  publisher = {Springer},
  year = {2000}
}

@book{milnor1963morse,
  author = {Milnor, John Willard},
  title = {{Morse} theory},
  number = {51},
  publisher = {Princeton University Press},
  year = {1963}
}

@book{matsumoto2002introduction,
  author = {Matsumoto, Yukio},
  title = {An introduction to {Morse} theory},
  volume = {208},
  publisher = {American Mathematical Soc.},
  year = {2002}
}

@book{kuznetsov1998elements,
  author = {Kuznetsov, Yuri A.},
  title = {Elements of Applied Bifurcation Theory},
  volume = {112},
  publisher = {Springer},
  year = {1998}
}

@article{keane2017climate,
  author = {Keane, Andrew and Krauskopf, Bernd and Postlethwaite, Claire M},
  title = {Climate models with delay differential equations},
  journal = {Chaos: An Interdisciplinary Journal of Nonlinear Science},
  volume = {27},
  number = {11},
  pages = {114309},
  year = {2017}
}

@book{shil2001methods,
  author = {Shil'nikov, Leonid P.},
  title = {Methods of qualitative theory in nonlinear dynamics},
  volume = {5},
  publisher = {World Scientific},
  year = {2001}
}

@book{wiggins2003introduction,
  author = {Wiggins, Stephen},
  title = {Introduction to Applied Nonlinear Dynamical Systems and Chaos},
  volume = {2},
  publisher = {Springer},
  year = {2003}
}

@book{audin2014morse,
  title     = {Morse Theory and Floer Homology},
  author    = {Audin, Mich{\`e}le and Damian, Mihai},
  year      = {2014},
  publisher = {Springer},
  note      = {Translated from the French by Reinie Ern{\'e}}
}

@article{terrien2023merging,
  author = {Terrien, Soizic and Krauskopf, Bernd and Broderick, Neil G. R. and Pammi, Venkata A. and Braive, R{\'e}my and Sagnes, Isabelle and Beaudoin, Gr{\'e}goire and Pantzas, Konstantinos and Barbay, Sylvain},
  title = {Merging and disconnecting resonance tongues in a pulsing excitable microlaser with delayed optical feedback},
  journal = {Chaos: An Interdisciplinary Journal of Nonlinear Science},
  volume = {33},
  number = {2},
  pages = {023142},
  year = {2023}
}

@article{peckham1995bananas,
  title={{B}ananas and banana splits: a parametric degeneracy in the {H}opf bifurcation for maps},
  author={Peckham, Bruce B. and Frouzakis, Christos E. and Kevrekidis, Ioannis G.},
  journal={SIAM Journal on Mathematical Analysis},
  volume={26},
  number={1},
  pages={190--217},
  year={1995},
  publisher={SIAM}
}

@article{peckham2002lighting,
  title={Lighting {A}rnold flames: resonance in doubly forced periodic oscillators},
  author={Peckham, Bruce B. and Kevrekidis, Ioannis G.},
  journal={Nonlinearity},
  volume={15},
  number={2},
  pages={405--428},
  year={2002},
  publisher={IOP Publishing}
}

@inproceedings{bolduc2026resonance,
  author    = {Bolduc-St-Aubin, Samuel and Krauskopf, Bernd},
  title     = {From Resonance to Chaos in a {DDE} Climate Model},
  booktitle = {24th IFAC World Congress, Busan, Republic of Korea},
  year      = {2026}
}

@article{bolduc2026enso,
  author  = {Bolduc-St-Aubin, Samuel and Subramanian, Priya and Krauskopf, Bernd},
  title   = {Resonance structure of a periodically forced delay differential equation model for the {El Ni\~{n}o}--{Southern Oscillation}},
  journal = {Preprint},
  year    = {2026},
  note    = {}
}
\end{document}